\documentclass[letter,11pt]{article}

\usepackage{amsmath}

\usepackage{comment}
\usepackage[linktocpage=true]{hyperref}
\usepackage[margin=1.2 in, top=1 in, bottom= 1.2 in]{geometry}

\usepackage[english]{babel}
\usepackage{float}
\usepackage{amsfonts}
\usepackage{mathrsfs}
\usepackage[mathscr]{euscript}
\usepackage{bbm}
\usepackage{latexsym}
\usepackage{amssymb}
\usepackage{mathtools}
\usepackage{graphicx}
\usepackage{tikz}
\usepackage{placeins}

\usepackage{enumitem}
\setlist{nolistsep}
\usepackage{amsthm}
\usepackage{relsize}
\usepackage{nicefrac}
\usepackage[capitalize]{cleveref}

\newtheoremstyle{plain}{3mm}{3mm}{\slshape}{}{\bfseries}{.}{.5em}{}
\newtheoremstyle{definition}{2mm}{2mm}{}{}{\bfseries}{.}{.5em}{}
\theoremstyle{plain}
	
\newtheorem{theorem}{Theorem}

\theoremstyle{definition}

\theoremstyle{plain}
\newtheorem*{namedthm}{\namedthmname}
\newcounter{namedthm}
	

\newcommand{\J}{\mathcal{J}}
\newcommand{\JJ}{\overset{\sim}{\mathcal{J}}}

\newcommand{\Z}{\mathbb{Z}}
\newcommand{\R}{\mathbb{R}}

\newcommand{\FF}{\textbf{F}}
\newcommand{\GG}{\textbf{G}}

\newcommand{\orbit}{$\text{Orb}_{F}(t)$}

\newcommand{\li}[1]{\overset\longrightarrow{#1}}

\title{Ergodicity in the linked twist maps with counter-oriented twists.}
\author{Aritro Pathak}

\begin{document}

\maketitle

\begin{abstract}
    We reduce the earlier known optimal twist parameter for which ergodicity is established in the linked twist map with two linear twists in opposite sense, in the most general setting. Further, here we obtain ergodicity with possibly only one-fold twists in either lobe, while earlier results only applied for twist parameters at least 2. Almost hyperbolicity is easily established for twist parameters greater than 2, while in the most general setting of the linked twist map with both twists of equal magnitude, ergodicity was earlier established in the most general setting by Przytycki for twist parameters greater than 4.15 with at least two-fold twists in each lobe. Here we reduce this optimal twist parameter to 3.47 in the general setting. These techniques can be effected to make further improvements when additional assumptions are made on the dimensions of the strips or when the linked twist map is modified in other natural ways.
\end{abstract}

\section{Introduction:}

\noindent The linked twist map is a classical and well studied dynamical system that exhibits hyperbolic behavior, like the Arnold's cat map or geodesic flow on a negative curvature surface; see \cite{Sp},\cite{Prz},\cite{BE},\cite{Woj}, \cite{Devaney} and the references therein. We significantly extend earlier techniques of Przytycki based on Pesin theory for singularities as modified by Katok and Strelcyn, to extend the range of twist parameters for which ergodicity is established for the classical linked twist map in a general setting, with twists in the two lobes that are counter oriented, while the twist parameters $k,m$ defined below are any non-zero integers but of opposite sign.  When the twists reinforce each other, the question has been dealt with by Burton and Easton \cite{BE}, as well as Wojtowski \cite{Woj}. Also see the Introduction of \cite{Liv} for a discussion of this situation. 

As we show below, hyperbolicity (or ``almost hyperbolicity" in the language of \cite{Prz}) is achieved for all twist parameters $\alpha >2$, and one typically expects ergodicity for the map to also be established for all $\alpha>2$. This has not yet been shown to be true. Upon markedly extending the methods of \cite{Prz} for this setting, we extend the set of twist parameters for which ergodicity is established. 
 
\begin{figure}
\centering
    \begin{tikzpicture}
         \draw[thick] (-2,0.5) to (2,0.5);
         \draw[thick] (-2,-0.5) to (2,-0.5);
         \draw[->,thick] (-2,0.5) to (-2,-0.5);
         \draw[->,thick] (2,0.5) to (2,-0.5);
         \draw[thick] (0.5,0.5) to (0.5,2);
         \draw[thick] (-0.5,0.5) to (-0.5,2);
         \draw[->,thick] (0.5,2) to (-0.5,2);
         \draw[thick] (0.5,-0.5) to (0.5,-2);
         \draw[thick] (-0.5,-0.5) to (-0.5,-2);
         \draw[->,thick] (0.5,-2) to (-0.5,-2);
         \draw[->] (-0.4,1.7) to (0.2,-1.9);
         \draw[->] (-1.9,-0.3) to (1.9,0.2);
         \draw[thick] (0.5,0.5) to (0.5,-0.5);
         \draw[thick] (-0.5,0.5) to (-0.5,-0.5);
         \node [above] at (-0.62,0.5) {$A$};
         \node [above] at (0.62,0.5) {$B$};
         \node [below] at (0.62,-0.5) {$C$};
         \node [below] at (-0.62,-0.5) {$D$};
         \node[right] at (2,0.5) {$P$};
         \node[right] at (2,-0.5) {$Q$};
         \node[above] at (0.5,2) {$R'$};
         \node[above] at (-0.5,2) {$S'$};
         \node[below] at (0.5,-2) {$R$};
         \node[below] at (-0.5,-2) {$S$};
         \node[left] at (-2,-0.5) {$Q'$};
         \node[left] at (-2,0.5) {$P'$};

         \draw[thick] (4,-1) to (4,1);
         \draw[thick] (6,-1) to (6,1);
         \draw[thick] (4,-1) to (6,-1);
         \draw[thick] (4,1) to (6,1);
         \node[right] at (6,1.1) {$B$};
         \node[left] at (4,1.1) {$A$};
         \node[right] at (6,-1.1) {$C$};
         \node[left] at (4,-1.1) {$D$};
         \draw[<-] (4.5,0.8) to (5.5,0.8);
         \draw[->] (4.5,-0.8) to (5.5,-0.8);
         \draw[->] (5.8,-0.5) to (5.8,0.5);
         \draw[<-] (4.2,-0.5) to (4.2,0.5);

    \end{tikzpicture}
\caption{The linked twist map shown in the figure. The horizontal track is the region $H=\{(0\leq x \leq 1)\times (y_0\leq y\leq y_1)\}$, the vertical track is the region $V=\{(x_0\leq x \leq x_1)\times (0\leq y\leq 1)\}$. The map is linked on the central square region $S$, by the twist map $F$ on $H$ and a vertical twist map $G$ on $V$. The points $A,B,C,D$ delimit the central square $S$. The pairs $(P,P'),(Q,Q'),(R,R'),(S,S')$ are identified and thus the segments $PQ,P'Q'$ are identified together and so are the segments $SR,S'R'$.The twist on the horizontal strip is denoted by $F$ and that on the vertical strip is denoted by $G$, although the underlying map is the same. The picture on the right shows the local twists near the boundary in the central square region $S$.}
\label{fig:figmain}
\end{figure}
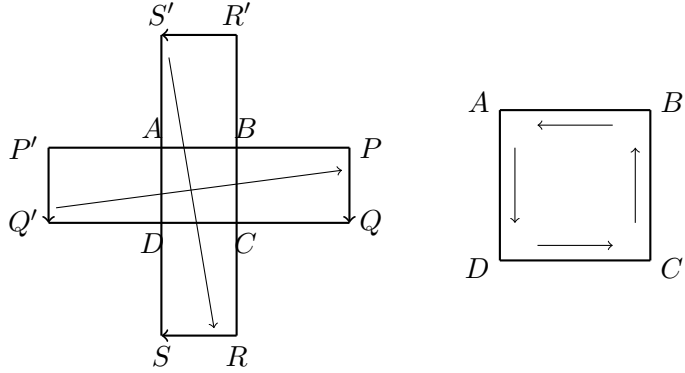

The methods established in \cite{Prz} using Pesin theory are canonical for the case where the two separate twisting tracks individually are homeomorphic to the set $[0,1]\times [0,1]/ \sim$ where each point $(x,0)$ is identified with $(x,1)$ for each $x\in [0,1]$. These two twisting tracks are then linked with each other. In other words, we have the case shown in the \cref{fig:figmain}. 

The dynamics can then be described as a product of two successive twists, one horizontally which we term as $F$, and the other vertical twist which we call $G$, so that the map becomes, with $\alpha>0,\beta>0$, in the domain $H=\{(x,y):x\in [0,1],y\in [y_0,y_1]\}$ for $F$ and the domain $V=\{(x,y): x\in [x_0,x_1],y\in [0,1] \}$ for $G$:

\begin{equation}\label{eq:eq2}
\FF\cdot
\begin{pmatrix}
x\\
y-y_0
\end{pmatrix}=\begin{pmatrix}
1 & \alpha \\
0 & 1 
\end{pmatrix} \cdot
\begin{pmatrix}
x\\
y-y_0
\end{pmatrix},
\end{equation}

\begin{equation}\label{eq:eq2}
\GG\cdot
\begin{pmatrix}
x-x_0\\
y
\end{pmatrix}=\begin{pmatrix}
1 & 0 \\
-\beta & 1 
\end{pmatrix} \cdot
\begin{pmatrix}
x-x_0\\
y
\end{pmatrix},
\end{equation}

The map $\FF$ is the identity when $y\in[0,1], y\notin [y_0,y_1]$ and $G$ is the identity when $x\in [0,1], x\notin [x_0,x_1]$. Further, we take the $\alpha(y_1 -y_0)=k$ and $\beta(x_1 -x_0)=m$ for some integers $k,m$. Following \cite{Prz}, we call the maps $\FF,\GG$ respectively $(k,\alpha)$ and $(m,\beta)$ twists.

We call the product $\Phi=\GG\circ \FF$. All the maps $\FF,\GG,\Phi$ preserve Lebesgue measure on $H\cup V$.

While we only deal with the linear case, the methods present here can approximate the case where:

\begin{align*}
   \FF(x,y)=(x+f(y),y ) \  \text{for} \ (x,y)\in H, \ \ \ \GG(x,y)= (x,y+g(x)), \ \text{for} \ (x,y)\in V,
\end{align*}

with $f:[y_0,y_1]\to \R$, $g:[x_0,x_1]\to \R$ both $C^{2}$ functions with $f(y_0)=g(x_0)=0$, $f(y_1)=k, g(x_1)=m$, with 

\begin{align*}
    \frac{df}{dy}\neq 0, \text{and} \ \ \frac{dg}{dx}\neq 0,
\end{align*}

for every $y\in [y_0,y_1]$ and $x\in [x_0,x_1]$. In this case one can define corresponding $\alpha$ and $\beta$ parameters as one of the extemum values taken by these derivatives of $f$ and $g$ on these closed intervals. For convenience, we only deal with the linear case, but the essential dynamics in the above case is also the same.

In ongoing work,\cite{Pat} a flow arising after a Dehn-type surgery on the unit tangent bundle of a genus two surface $\mathbb{S}M$ is constructed, where the original flow prior to surgery consists of a simple periodic flow along each of the fibers of the unit tangent bundle, and upon a Dehn surgery along an annulus around a curve $C$ on $\mathbb{S}M$ where $C$ projects to a closed self intersecting geodesic on $M$, we can project the flow to a linked twist map where the identification of the boundaries is such that the right edge of the horizontal track is identified with the top edge of the vertical track, and the bottom edge of the vertical track is identified with the left edge of the horizontal track.Here we have twists that are uniform $C^{0}$ at the boundary, and we have to modify the arguments of the original setting of \cite{Prz} because of this altered boundary identification.\footnote{See \cite{FH2013},\cite{FH2021} for some background information on this question.}

In this paper, we deal only with the well-studied canonical boundary identifications of the linked twist map as in \cref{fig:figmain}, and improve the value of the critical paramter $\alpha_0$ so that ergodicity is established, in this most general setting, for all $\alpha>\alpha_0$. \footnote{A similar analysis could be done for the case of maps linked several times, or where the identifications of the boundary are as in the Dehn-surgered problem in the upcoming manuscript \cite{Pat}, by adopting and extending the same idea used here.}

\begin{theorem}
    Consider a linked twist map $\Phi=\GG\circ \FF$ composed of a $(k,\alpha)$ horizontal twist denoted by $F$ and a $(m,\beta)$ vertical twist denoted by $G$, where $k$ and $m$ are any non-zero integers which have opposite signs, and $\alpha\beta\geq 12.04$, then $\Phi$ is ergodic.
\end{theorem}

Henceforth, we deal with the case where $\alpha=\beta$, without loss of generality, and our critical parameter $\alpha_0= 3.47$.

Our method establishes ergodicity; we are not able to establish the Bernoulli property for this reduced critical parameter. However, we are able to remove the restriction in the hypothesis of \cite{Prz} that $|k|,|m|\geq 2$.

For the case of twists having opposite signs, and where $\alpha\beta>4$, we get hyperbolicity for the map $\Phi$.

\begin{align*}
    D\Phi= D\GG\circ D\FF=\begin{pmatrix}
                        1 & 0 \\
                        -\beta & 1
                        \end{pmatrix}\circ \begin{pmatrix}
                        1 & \alpha \\
                        0 & 1
                       \end{pmatrix}=\begin{pmatrix}
                        1 & \alpha \\
                        -\beta & 1-\alpha\beta
                       \end{pmatrix}
\end{align*}

In this case, the eigenvalues are:

\begin{align*}
    \frac{(2-\alpha\beta)\pm\sqrt{(\alpha\beta)^{2}-4\alpha\beta}}{2},
\end{align*}

which shows hyperbolicity is achieved for $\alpha\beta>4$.

Without loss of generality, we can take the twists in the two lobes to be of equal magnitude, i.e. $|\alpha|=|\beta|$, since otherwise, we can rescale one of the two variables by a factor of $\sqrt{\alpha/\beta}$ so that in the subsequent analysis, the twists in both the lobes are of magnitude $\sqrt{\alpha\beta}$ and of opposing signs. \footnote{While effecting this change of variables, the product $\alpha\beta$ remains constant and so does the eigenvalues and eigenvectors above.}

In case both the twist parameters are of the same magnitude $|\alpha|=|\beta|$, the statement of the theorem reduces to saying that the above critical parameter is $\sqrt{12.04}=3.47$. Henceforth we make this assumption of the two twists having the same magnitude.

As noted in \cite{Prz}, the expanding eigenvector $(\xi_1,\xi_2)$ satisfies:

\begin{align*}
    \frac{\xi_1}{\xi_2}=-\Big( \frac{\alpha}{2} \Big) +\sqrt{(\frac{\alpha}{2})^{2}-1},
\end{align*}

and henceforth, we call this number $L_\alpha$. As explained in the next section, the line in $\R^{2}$ with slope $1/L_\alpha$, and the vertical line, together delimit the vertical cone under consideration.

\bigskip

In several places it can be seen why the critical value for the optimal twist of $3.47$ is around the best possible with our refinements. We refer the reader to Figure 5. Here it would be enough to require $|\li{JA}|$ to be around $(2/3)l_{v}(\gamma)$ and also $\alpha \gtrsim 3.3$: in this case $|\li{BE}|$ is either at least $\l_{v}(\gamma)$ and we are done, and otherwise $|\li{EF}|=|L_\alpha| |\li{BE}|\lesssim (1/3)l_{v}(\gamma)$ and then $|\li{FP}|\geq |\li{JA}|-|\li{EF}|\gtrsim (1/3)l_{v}(\gamma)$. In this case we have the possibility of having $|\li{GK}|\approx (1/3)l_{v}(\gamma)(3.3- 0.3)\approx l_{v}(\gamma)$ or the segment gets cut off by the upper edge of $S$, and our argument continues. The same argument holds for the segments that enter $S$ while being cut-off by the right edge of $S$ as well and further on similar considerations apply in the analysis for Figure 8.


\section{Proof of Theorem 1:}

We use Pesin theory for maps having singularities. For reference, check the Appendix of \cite{Prz}, and \cite{Kat}. Note that Lebesgue almost every point of $S$ returns to $S$ with positive frequency. We need to verify for all such points $x,y\in S$, given the local unstable manifold $\gamma^{u}(x)$ through $x$ and the local stable manifold $\gamma^{s}(y)$ through $y$, that there exist integers $m,n$ large enough such that $\Phi^{m}(\gamma^{u}(x))\cap \Phi^{-n}(\gamma^{s}(y))\neq 0$.

 It is enough to work in Subsections 2.1,2.2,2.3, on the dynamics of the forward iterate of $\gamma^{u}(x)$ for almost every $x$ in the square, and the corresponding analysis for the backward iterates of $\gamma^{s}(y)$ is similar. 

\begin{figure}
\centering
    \begin{tikzpicture}
         \draw[thick] (1,-1) to (1,1);
         \draw[thick] (3,-1) to (3,1);
         \draw[thick] (1,-1) to (3,-1);
         \draw[thick] (1,1) to (3,1);
         \draw[thick] (1,0) to (1.5,-1);
         \draw (3,-0.3) to (1.5,-1);
         \draw[thick] (3,-0.3) to (2,1);
         \draw (1,0.3) to (2,1);
         \draw[thick] (1,0.3) to (1.3,-1);
         \draw (3,-0.6) to (1.3,-1);
         \draw[thick] (2.6,1) to (3,-0.6);
         \draw (2.6,1) to (1,0.45);
         \draw (2.65,1) to (1,0.4);
         \draw (1.3,-1) to (3,-0.65);
         \draw[thick] (1,0.45) to (1.35,-1);
         \draw[thick] (2.65,1) to (3,-0.65);
         \draw (2.65,1) to (1,0.45);
         \draw (1.35,-1) to (3,-0.70);
         
         \node[right] at (3,1.1) {$B$};
         \node[left] at (1,1.1) {$A$};
         \node[right] at (3,-1.1) {$C$};
         \node[left] at (1,-1.1) {$D$};
         \node[above] at (1.5,-0.8) {$\gamma_{1_1}$};
         \node[above] at (2.5,-0.3) {$\gamma_{2_1}$};
         \node[above] at (2.3,-0.7) {$\gamma_{1_2}$};
         \node[above] at (2,-2.5) {$(a)$};
         \node[left] at (1,0) {$P_1$};
         \node[below] at (1.5,-1) {$P_2$};
         \node[right] at (3,-0.3) {$P_3$};
         \node[above] at (2,1) {$P_4$};

         \draw[thick] (10,-1) to (10,1);
         \draw[thick] (12,-1) to (12,1);
         \draw[thick] (10,-1) to (12,-1);
         \draw[thick] (10,1) to (12,1);
         \draw[thick] (11,-1) to (10.70,1);
         \node[above] at (10.7,-0.4) {$\gamma_1$};
         \node[right] at (12,1.1) {$B$};
         \node[left] at (10,1.1) {$A$};
         \node[right] at (12,-1.1) {$C$};
         \node[left] at (10,-1.1) {$D$};
         \node[above] at (11,-2.5) {$(c)$};

         \draw[thick] (6,-1) to (6,1);
         \draw[thick] (8,-1) to (8,1);
         \draw[thick] (6,-1) to (8,-1);
         \draw[thick] (6,1) to (8,1);
         \draw (6,-0.5) to (8,0);
         \draw[thick] (6,-0.5) to (6.15,-1);
         \draw[thick] (8,0) to (7.75,1);
         \draw[thick] (6.9,1) to (7.3,-1);
         \node[right] at (8,1.1) {$B$};
         \node[left] at (6,1.1) {$A$};
         \node[right] at (8,-1.1) {$C$};
         \node[left] at (6,-1.1) {$D$};
         \node[above] at (6.5,-0.4) {$\gamma_1$}; 
         \node[above] at (7,-2.5) {$(b)$};
       \node[above] at (7.5,0.5) {$\gamma_2$}; 
       \node[above] at (6.4,-1) {$\gamma_3$}; 
       \node[above] at (7.4,-1) {$\gamma_4$};

    \end{tikzpicture}
\caption{The dynamics for the forward iterates of the unstable manifold. In part(a), for each $i\in \{1,2,\dots,\}$, we have a sequence of segments $\gamma_{i_1},\gamma_{i_2} \subset \Phi^{m_i}(\gamma^{u}(x)),\FF\circ \Phi^{m_i}(\gamma^{u}(x))$ respectively, each of which intersect two adjacent edges of $S$; the segments $\gamma_{1_2}, \gamma_{2_1}$ belong to $\FF\circ\Phi^{m_0 +1}(\gamma^{u}(x)), \Phi^{m_{0}+2}(\gamma^{u}(x))$ and so on. In part (b), we have a horizontal segment $\gamma_1 \subset \FF\circ \Phi^{m_0}(\gamma^{u}(x))$ through $S$ and then subsequent vertical segments $\gamma_2,\gamma_3,\gamma_4 \subset \Phi^{m_0 +1}(\gamma^{u}(x))$and in part (c), we have a vertical segment $\gamma_1$ belonging to $\Phi^{m_0}(\gamma^{u}(x))$ through $S$, for some positive integer $m_0$, in which case we are done.  }
\label{fig:figsec}
\end{figure}
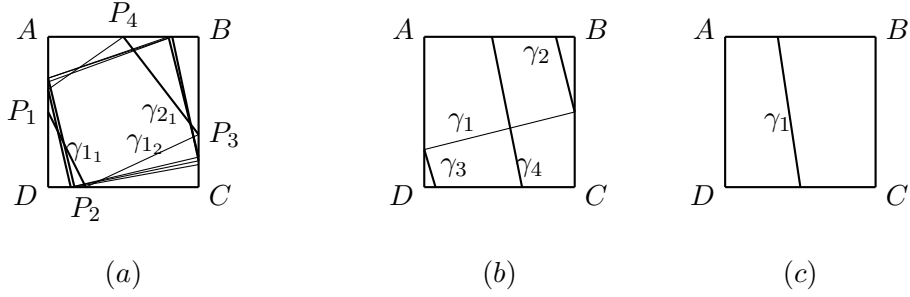

The outcome is that, for any $x\in S$, either of the following two cases happen:
\begin{enumerate}

    \item Either a horizontal segment through $S$ in some iterate $F\circ \Phi^{m_0}(\gamma^{u}(x))$, or a vertical segment through $S$ in some iterate $\Phi^{m_0}$ and these cases are shown in Figure 2(b) and 2(c).

    \item A contiguous union of a countable infinity of segments, each of which touch two adjacent sides as shown in \cref{fig:figsec}(a). We enumerate these segments as $\gamma_{i_1}, \gamma_{i_2}$ with $i\in \{1,\dots,\infty\}$, and such that $\gamma_{i_2}\subset \FF(\gamma_{i_1}), \gamma_{(i+1)_1}\subset \GG(\gamma_{i_2})$ for each $i\in \{1,\dots, \infty \}$.
\end{enumerate}

There are two corresponding cases for the backward iterates of the unstable manifold.

In subsection (2.4), we combine the above two situations coming from each of the cases of the forward iterates of the unstable manifold along with the backward iterates of the stable manifold, to get the following four possibilities:

\begin{enumerate}
    \item[a.] Situation 1 for the forward iterate of the unstable manifold $\gamma^{u}(x)$, and also situation 1 for the backward iterate of the stable manifold $\gamma^{s}(y)$, which means either a horizontal segment through $S$ belonging to $\FF\circ \Phi^{m_0}(\gamma^{u}(x))$, or a vertical segment through $S$ belonging to $\Phi^{m_0}(\gamma^{u}(x))$, and either a horizontal segment through $S$ belonging to $ \Phi^{-n_0}(\gamma^{s}(y))$ or a vertical segment belonging to $\GG^{-1}\circ \Phi^{-n_0}(\gamma^{s}(y))$.

    \item[b.] Situation 2 for the forward iterates of the unstable manifold and situation 1 for the backward iterates of the stable manifold.

    \item[c.] Situation 1 for the forward iterates of the unstable manifold and situation 2 for the backward iterates of the stable manifold.

    \item[d] Situation 2 for both the forward iterate of the unstable manifold and backward iterates of the stable manifold.
    
\end{enumerate}

\bigskip

In each of these four contingencies above, we will show that above the critical twist parameter, there is a vertical segment belonging to $\Phi^{m}(\gamma^{u}(x))$ for some large enough $m$ that intersects a horizontal segment belonging to $\Phi^{-n}(\gamma^{s}(y))$ for some large enough $n$, thus establishing ergodicity for the linked twist map above the critical twist parameter. 

\bigskip


 \bigskip

\subsection{Initial argument}

We reproduce the argument of \cite{Prz} with a more detailed exposition, with reference to the two diagrams in \cref{fig:fifth}, and start our argument after this. Following \cite{Prz}, we denote the first return to $S$ under the maps $\FF,\GG,\Phi$ respectively as $\FF_s, \GG_s, \Phi_s$, and then fixing any segment $\gamma\in \Phi_{s}^{r}(\gamma^{u}(x))$ in the $r$'th iterate of the first return map $\Phi_{s}$, which lies within the specified cone, it is enough to focus on the first return to $S$ under the horizontal twist map of an initial segment $\gamma\in S$ that lies within its permissible cone, i.e. on $\FF_{s}(\gamma)$  the first time $m_1\in \Z_{+}$ that $\FF^{m_1}(\gamma)\cap S \neq 0$. 

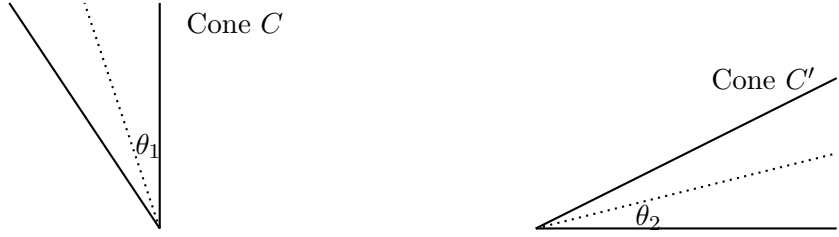
\begin{figure}
\centering
    \begin{tikzpicture}
        \draw[thick] (-7,0) to (-9,3);
        \draw[dotted, thick] (-7,0) to (-8,3);
        \draw[thick] (-7,0) to (-7,3);
        \node [below] at (-6,3) {Cone $C$};
        \draw[thick] (-2,0) to (2,2);
        \draw[dotted, thick] (-2,0) to (2,1);
        \draw[thick] (-2,0) to (2,0);
        \node [right] at (0.2,2) {Cone $C'$};
        \node [above] at (-7.15,0.8) { $\theta_1$};
        \node [above] at (-0.5,-0.15) { $\theta_2$};
        
    \end{tikzpicture} 
    \caption{The slope of the dotted segment in the vertical Cone $C$ in the left is taken as $-\tan \theta_1$ and the slope of the dotted segment in the horizontal Cone $C'$ on the right is taken as $-\tan \theta_2$. Whenever we speak of the slope of the segment in the cone in which the segment lies, we understand from context whether the cone is horizontal or vertical. We have uniform bounds of $L_\alpha\leq L\leq 0$ for the slope $L$ of the segment within the cone, with $L_{\alpha}=-\frac{\alpha}{2}+\sqrt{(\frac{\alpha}{2})^{2}-1}$.}
    \label{fig:fourth}
\end{figure}

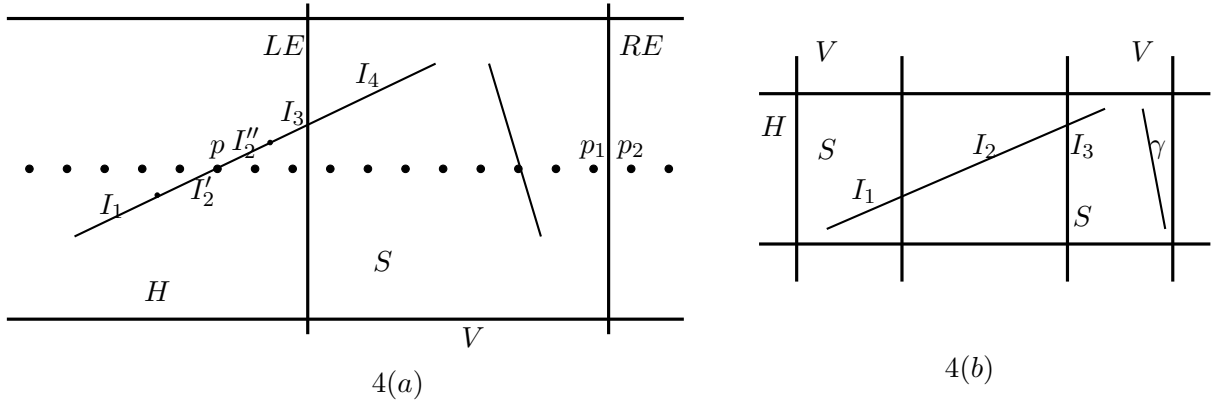
\begin{figure}
\centering
    \begin{tikzpicture}
        \draw[very thick] (2,2) to (-7,2);
        \draw[very thick] (2,-2) to (-7,-2);
        \draw[very thick] (-3,-2.2) to (-3,2.2);
        \draw[very thick] (1,-2.2) to (1,2.2);
        \node [above] at (-2,-1.5) {$S$};
        \node [above] at (-5,-1.9) {$H$};
        \node [above] at (-0.8,-2.5) {$V$};
        \node [above] at (-1.8,-3.2) {$4(a)$};
        \node [above] at (5.8,-3.0) {$4(b)$};
        \node [above] at (-4.2,0) {$p$};
        \draw [thick] (-6.1,-0.9) to (-1.3,1.4);
        \draw [thick] (-0.59,1.4) to (0.1,-0.9);
        \draw [fill] (-6.7,0) circle [radius=0.05];
        \draw [fill] (-6.2,0) circle [radius=0.05];
        \draw [fill] (-5.7,0) circle [radius=0.05];
        \draw [fill] (-5.2,0) circle [radius=0.05];
        \draw [fill] (-4.2,0) circle [radius=0.05];
        \draw [fill] (-4.7,0) circle [radius=0.05];
        \draw [fill] (-3.7,0) circle [radius=0.05];
        \draw [fill] (-3.2,0) circle [radius=0.05];
        \draw [fill] (-2.7,0) circle [radius=0.05];
        \draw [fill] (-2.2,0) circle [radius=0.05];
        \draw [fill] (-1.7,0) circle [radius=0.05];
        \draw [fill] (-1.2,0) circle [radius=0.05];
        \draw [fill] (-0.7,0) circle [radius=0.05];
        \draw [fill] (-0.2,0) circle [radius=0.05];
        \draw [fill] (0.3,0) circle [radius=0.05];
        \draw [fill] (0.8,0) circle [radius=0.05];
        \draw [fill] (1.3,0) circle [radius=0.05];
        \draw [fill] (1.8,0) circle [radius=0.05];
        \node [left] at (-2.9,1.65) {$LE$};
        \node [right] at (1,1.65) {$RE$};
        \node [above] at (0.8,0) {$p_1$};
        \node [above] at (1.3,0) {$p_2$};
        \node [left] at (-5.3,-0.5) {$I_{1}$};
        \draw[fill] (-5.0,-0.35) circle [radius=0.03];
        \draw[fill] (-3.5,0.35) circle [radius=0.03];
        \node [right] at (-4.7,-0.3) {$I'_{2}$};
        \node [left] at (-3.5,0.35) {$I''_{2}$};
        \node [above] at (-3.2,0.45) {$I_{3}$};
        \node [above] at (-2.2,0.95) {$I_{4}$};

        \draw[very thick] (3,1) to (9,1);
        \draw[very thick] (3,-1) to (9,-1);
        \draw[very thick] (3.5,-1.5) to (3.5,1.5);
        \draw[very thick] (4.9,-1.5) to (4.9,1.5);
        \draw[very thick] (7.1,-1.5) to (7.1,1.5);
        \draw[very thick] (8.5,-1.5) to (8.5,1.5);
        \draw[thick] (3.9,-0.8) to (7.6,0.8);
        \draw[thick] (8.4,-0.8) to (8.1,0.8);
        \node [above] at (8.3,0) {$\gamma$};
        \node [above] at (3.9,0) {$S$};
        \node [above] at (7.3,-0.9) {$S$};
        \node [above] at (3.2,0.3) {$H$};
        \node [above] at (3.9,1.3) {$V$};
        \node [above] at (8.1,1.3) {$V$};
        \node [left] at (4.7,-0.3) {$I_1$};
        \node [below] at (6.0,0.6) {$I_2$};
        \node [below] at (7.3,0.6) {$I_3$};

    \end{tikzpicture} 
    \caption{Reproduction of the diagrams of \cite{Prz} relevant to the dynamics in our problem. We denote the union of the segments $I'_2 \cup I''_2=I_2$. These same segments are depicted again in \cref{fig:first} below and labelled, and those labelings are used in the paper throughout.}
    \label{fig:fifth}
\end{figure}

 We also use the terminology of identifying the slope of a segment within the cone in which it lies; this is illustrated and explained in \cref{fig:fourth}.
 
We have the following four cases:

\begin{enumerate}
    \item[(i)] $\FF^{m_{1}}(\gamma) $ contains a horizontal segment. The analysis here is subsumed in the analysis done for Case (ii) below.

    \item[(ii)] The right side of $\FF^{m_{1}}(\gamma) $ intersects $S$ but we do not have a horizontal segment. This case is shown in \cref{fig:fifth}(a). 

    \item[(iii)] The left side of $\FF^{m_{1}}(\gamma) $ intersects $S$ but we do not have a horizontal segment. This case is entirely analogous to the above one, and we will only discuss the above case.

    \item[(iv)] Both sides of $\FF^{m_{1}}(\gamma) $ intersect $S$, which is depicted in \cref{fig:fifth}(b).

\end{enumerate}

\bigskip

 We first deal with the Case (ii) above. As what turns out to be the canonical approach as in \cite{Prz}, we divide $\FF^{m_1}(\gamma)\setminus S$ into three intervals $I_1,I_2,I_3$ which will be described later. Also, we denote $\FF^{m_1}\cap S=I_4$.
 
 We identify a point $p\in I_2$ that has a periodic orbit under $F$ with distance between points of orbit being $d$. We define the distance $d$ as follows with reference to the segments shown in \cref{fig:fifth}: look for the unique integer $q$ such that $1/q < \alpha l_v (I_2)$ and $1/(q-1) \geq \alpha l_v (I_2)$, which is $q= \lfloor \frac{1}{\alpha l_v (I_2)} +1 \rfloor$. Given the segment $I_2$, let the vertical endpoints of $I_2$ be $y_1,\text{and} \ y_2=y_1+l_{v}(I_2)$. Under a forward iterate of the map, one endpoint moves forward by a distance of $\alpha y_1$ and the other endpoint moves forward by a distance $\alpha(y_1 +l_v(I_2))$. We seek a point on $I_2$ such that this point moves under this iteration of the map by a rational amount of $t/q$ with $\alpha y_1 \leq t/q < \alpha(y_1 +l_v(I_2))$, with $q=\lfloor \frac{1}{\alpha l_v (I_2)} +1 \rfloor$ as above. Clearly such a point would exist. The period of such a orbit, depending on whether we find $t,q$ coprime, is some divisor of $q$ and the distance between succesive points of the orbit of $t$ is some multiple of $d=\frac{1}{\lfloor \frac{1}{\alpha l_v (I_2)} +1 \rfloor}$ and thus at least this value. While we could have chosen a larger value of $q$, that would have made the distance $d$ smaller, which would need a stronger twist. Henceforth the orbit of $t$ under the map $\FF$ is denoted by \orbit.

The last point of \orbit \ just to the left of the right edge $RE$ of $S$, as shown in Fig 4(a) is denoted by $p_1$, and the point just to the right of $RE$ is denoted by $p_2$. We denote by $m_2>0$ the first time the when $F^{m_2}(p)$ lies between $p$ and $RE$. This cannot include the point $p$ since then by definition \orbit \ is disjoint from $S$, but which is a contradiction since $F^{-m_{1}}\in \gamma \subset S$. 

Now consider the segment $\mathcal{J}_0=I^{''}\cup I_3$, and then further on, $\mathcal{J}_m=\FF(\mathcal{J}_{m-1}\setminus S)$ for all $m=1,\dots,m_2$, i.e. at each stage up to the $m_2$'th, we consider the forward iterate of the previous segment after the portion of this previous segment already intersecting $S$ is excised. 

Then for each $m=1,2,\dots,m_2$, we have 

\begin{align}\label{eqverify}
    l_h(\J_m)\geq \text{min}(d+l_h(I_{2}'' \cup I_3),l_h(I_{2}'' \cup I_3)+\alpha l_v(I_{2}'' \cup I_3)).
\end{align}

For any fixed $1\leq m< m_2$, first term on the right comes from the case where $F^{w+1}(p)$ is actually the point of \orbit \  just to the left of $p$ and where the segment $\FF(\J_m)\cap S\neq 0$. In that case, $\J_m\setminus S$ is cut off by $LE$ and has length exactly $(d+l_h(I_{2}'' \cup I_3)$, and $\J_{m+1}=\FF(\J_m\setminus S)$ is larger than this length.\footnote{It will be arbitrarily larger depending on the strength of the twist $\alpha$.} On the other hand, when $\J_{m}\cap S=\phi$, then we get a segment of horizontal length at least $l_h(I_{2}'' \cup I_3)+\alpha l_v(I_{2}'' \cup I_3)$ in $\J_{m+1}$\footnote{This length could be significantly longer if we have a set of indices $J$ where the segments $\J_w$ do not intersect $S$, for any $w$ belonging to $J$.}.

When $\FF^{m_2}(p)$ lies between $p$ and $LE$, we must have a horizontal length at least $d$ being inserted into $S$ in the $m_2$'th iterate, and so $l_h(\J_{m_{2}}\cap S)\geq d$. That this is possible is easily verified from \cref{eqverify}. In fact there is going to be a horizontal inserted length bigger than $d$ which we cannot quantify accurately, and so we work with a lower bound of $d$. Furthermore, in this case this inserted segment $J_{m_{2}}\cap S$ touches the left edge $LE$. Further, if $\FF^{m_2}(p)=S\setminus\{p_1\}$, then we have:

\begin{align}\label{eqanother}
    l_h(\J_{m_2}\cap S)\geq \text{max} (d,l_h(I_{2}'' \cup I_3)+\alpha l_v(I_{2}'' \cup I_3))
\end{align}

Here the first term comes from the contingency where $\FF^{m_2}(p)$ is a point in $S$ just to the left of $p_1$ and $p_1$ is arbitrarily close ot $RE$, and $\J_{m_2}\cap S$ touches $RE$. Otherwise there is a segment of length $l_h(I_{2}'' \cup I_3)+\alpha l_v(I_{2}'' \cup I_3)$ within $S$, not necessarily touching the side $RE$.

Now assume that $\FF^{m_2}(p)=p_1$, and following \cite{Prz}, denote the distance of $p_1$ to $RE$ by $\tau\cdot d$. So if \cref{eqanother} is not satisfied, we have $l_h(\J_{m_2}\cap S)=\tau\cdot d$ and $\J_{m_2}\cap S$ touches $RS$.

Define $\JJ_{0}=I_1 \cup I_{2}'$ and then $\JJ_m =\FF(\JJ_{m-1}\setminus S)$ for $m=1,2,\dots m_2$. In this case, we have either 

\begin{align}\label{eqanotherr}
    l_{h}(\JJ_{m_2}\cap S)\geq \text{max}((1-\tau)\cdot d,l_h(I_{1} \cup I_2')+\alpha l_v(I_{1}\cup I_2')),
\end{align}

or $\JJ_{m_2}\cap S$ touches $LS$ with its left end, where the number $(1-\tau)\cdot d$ is used since we might have a case where for some $0<m<m_2$, we have $F^{m}(p)=p_2$ and the segment on the outside of $S$ gets cut off by $LE$ at this $m$'th step.\footnote{Again we actually have a longer segment than this within $S$ which we can't quantify, and so work with a lower bound of $(1-\tau)\cdot d$.} Since in this Case(ii), we don't assume a horizontal segment, from the above two contingencies, we can have either one of:

\begin{enumerate}
    \item $\JJ_{m_2}\cap S$ touches $LE$ with its left end and also \cref{eqanother} is true.

    \item $l_h(\JJ_{m_2}\cap S)=\tau\cdot d$ and $\J_{m_2}\cap S$ touches $RE$ and also \cref{eqanotherr} is true.
\end{enumerate}

In either case, whenever the optimum length inserted inside $S$ is lower bounded by $d$, this length necessarily touches either of the $LE$ or $RE$.

Note also that $\FF^{m_1}(\gamma)\cap S$ which we denote by $I_4$ touches the $LE$ as well.

\bigskip

The improvements in this paper concern with these above two contingencies where we have segments entering $S$ while necessarily touching either of the $LE$ or $RE$. We show that these segments that enter $S$ touching either of the left or the right sides, do not need to have a horizontal projection length as long as $\delta l_{v}(\gamma)$ for some $\delta>1$ arbitrarily close to $1$ as was needed in \cite{Prz}.\footnote{See Eq (3),(4),(5) on Page 350 in \cite{Prz}},
and can instead have lengths some $\beta_{i} l_{v}(\gamma), (i=1,2)$, with $\beta_i$ large enough to be determined later, but still $\beta_i<1, (i=1,2)$, and upon iterating successively with the maps $\FF,\GG$ further, we eventually get a segment within $S$ long enough of length $\delta l_{v}(\gamma)$ with now $\delta>1$, or we get a contiguous set of segments that touches all four edges, as shown in \cref{fig:figsec}(a,b)\footnote{As will be clear in the ensuing argument, we wish to ensure that either the case shown in \cref{fig:figsec}(a) appears, while otherwise we  get a complete horizontal or vertical segment in one more iteration of the map $\Phi$, which is the situation in the case shown in \cref{fig:figsec}(b).}. In the second case we are done, while in the first case, we keep on iterating the procedure to get longer segments that enter $S$ as each stage, until we also reach the second case, in which case we are done. 

It would thus be enough to require all of:

\begin{align}\label{deltareference}
    d\geq \delta \beta_2 \cdot l_{v}(\gamma) \\
    \alpha \cdot l_v(I_3) +l_h(I_3)\geq \delta \cdot l_{v}(\gamma)\\
    \alpha \cdot l_v(I_I) +l_h(I_I)\geq \delta \cdot l_{v}(\gamma)
\end{align}

or we would be also done if: 

\begin{align}
    l_{h}(I_4)\geq \delta \beta_1 \cdot l_{v}(\gamma).
\end{align}

\subsection{Part 1}

Below, we find optimum lower bounds on the parameters $\beta_1,\beta_2$, and use them to find the optimal parameter $\alpha$. Later, we work with Case(iv) and find the optimal value coming from the constraints there, which will be seen to be lower than the one coming from Case(ii) in the general case.

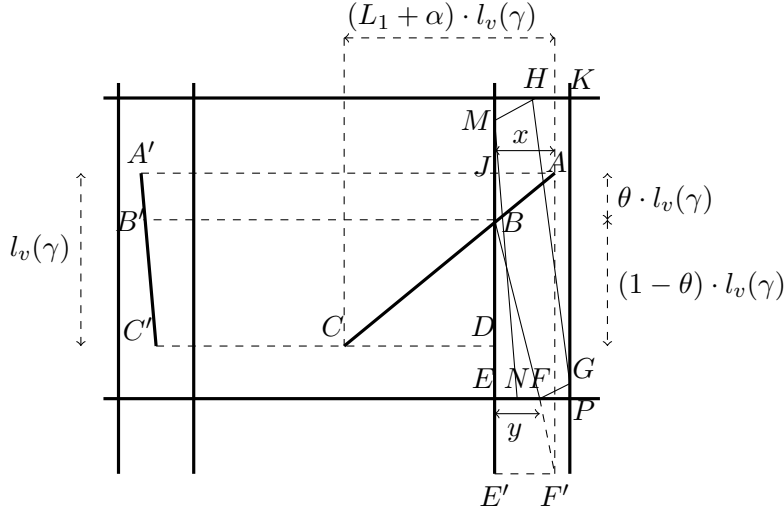
\begin{figure}
\centering
    \begin{tikzpicture}
        \draw[very thick] (2,2.2) to (2,-3);
        \draw[very thick] (1,2.2) to (1,-3);
        \draw[very thick] (-3,2.2) to (-3,-3);
        \draw[very thick] (-4,2.2) to (-4,-3);
        \draw[very thick] (-4.2,2) to (2.4,2);
        \draw[very thick] (-4.2,-2) to (2.4,-2);
        \draw[very thick] (-3.5,-1.3) to (-3.7,1);
        \draw[very thick] (-1,-1.3) to (1.8,1);
        \draw (1,0.38) to (1.6,-2);
        \draw[dashed] (1.6,-2) to (1.8,-3);
        \draw (1.6,-2) to (2,-1.8);
        \draw (2,-1.8) to (1.5,2);
        \draw[dashed] (-3.5,-1.3) to (1,-1.3);
        \draw[dashed] (1,-3) to (1.8,-3);
        \draw[dashed] (-3.7,1) to (1.8,1);
        \draw[dashed] (1,0.38) to (-3.67,0.38);
        \draw[<->] (1,1.3) to (1.8,1.3);
        \draw[<->] (1,-2.2) to (1.6,-2.2);
        \node [below] at (1.35,1.70) {$x$};
        \node [below] at (1.27,-2.2) {$y$};
        \node [above] at (1.82,0.88) {$A$};
        \node [right] at (0.95,0.35) {$B$};
        \node [above] at (-1.15,-1.3) {$C$};
        \node [above] at (0.85,-1.3) {$D$};
        \node [above] at (0.85,-2) {$E$};
        \node [above] at (1.6,-2) {$F$};
        \node [right] at (1.9,-1.6) {$G$};
        \node [right] at (1.9,-2.15) {$P$};
        \node [above] at (1.55,2) {$H$};
        \draw (1,1.7) to (1.55,2);
        \node [left] at (1.1,1.7) {$M$};
        \draw (1,1.7) to (1.3,-2);
        \node [above] at (1.3,-2) {$N$};
        \node [left] at (1.1,1.1) {$J$};
        \node [above] at (-3.7,1) {$A'$};
        \node [left] at (-3.5,0.38) {$B'$};
        \node [above] at (-3.74,-1.35) {$C'$};
        \node [below] at (1,-3) {$E'$};
        \node [below] at (1.8,-3) {$F'$};
        \draw [dashed] (1.8,-3) to (1.8,1);
        \draw [dashed,<->] (-4.5,-1.3) to (-4.5,1);
        \node [left] at (-4.5,0) {$l_{v}(\gamma)$};
        \node [right] at (1.85,2.25) {$K$};
        \node [right] at (2.5,0.65) {$\theta \cdot l_{v}(\gamma)$};
        \draw [dashed,<->] (2.5,0.38) to (2.5,1);
        \draw [dashed,<->] (2.5,0.38) to (2.5,-1.3);
        \node [right] at (2.5,-0.5) {$(1-\theta) \cdot l_{v}(\gamma)$};
        \draw [dashed,<->] (-1,2.8) to (1.8,2.8);
        \draw [dashed] (-1,2.8) to (-1,-1.3);
        \draw [dashed] (1.8,2.8) to (1.8,0.88);
        \node [above] at (0.3,2.75) {$(L_{1}+\alpha)\cdot l_{v}(\gamma)$};

    \end{tikzpicture} 
    \caption{Let $\gamma=\overset\longrightarrow{A'C'}$ be the initial segment within the square $S$. Here $m_1$ is the least positive integer such that $\FF^{m_1}(\gamma) \cap S \neq 0$ again, and $\FF^{m_1}(\gamma)=\overset\longrightarrow{CA}$. For the sake of clarity, we lift to $\R$ and this segment is depicted as intersecting $S$ on the right, corresponding to Case(ii).  Note that $|L_1|\leq |L_{\alpha}|=|-\frac{\alpha}{2}+\sqrt{(\frac{\alpha}{2})^{2}-1}|$. Also note that $\overset\longrightarrow{JA}=\beta_1 l_{v}(\gamma)$.}
    \label{fig:first}
\end{figure}

\bigskip

Below, in \cref{fig:first}, we redraw \cref{fig:fifth}(a) according to our convenience, and where the segment $\gamma$ is shown to be in $S$ but in a separate component in the lift of $H$ to $\R^{2}$. We have that $\overset\longrightarrow{AB}=I_4$ and also $\overset\longrightarrow{BC}=I_1\cup I_2 \cup I_3$.

We see from \cref{fig:first} that $y=|\overset\longrightarrow{EF}|$, $|\overset\longrightarrow{JB}|=\theta l_{v}(\gamma), |\overset\longrightarrow{BD}|=(1-\theta)l_{v}(\gamma)$, $l_{h}(\FF^{m}(\gamma))=|\overset\longrightarrow{JA}|+|\overset\longrightarrow{CD}|=(L_{1}+\alpha)l_{v}(\gamma)$, $x=|\overset\longrightarrow{JA}|=|\overset\longrightarrow{E'F'}|$, and with elementary geometry that:

\begin{align}
x= |\overset\longrightarrow{JA}|=\theta l_{v}(\gamma) (L_{1}+\alpha)=\beta_1 l_{v}(\gamma)\\
y(L_{2}+\alpha)=|\overset\longrightarrow{BE}|
\end{align}

Here $\theta$ is some number between $0$ and $1$, we defined $\beta_1=\theta(L_{1}+\alpha)$, $L_{1}$ defined in the figure is the slope of $\overset\longrightarrow{A'C'}$ within the cone $C$ in which it lies. Further, $L_{2}$\footnote{Here too, as usual, we have $L_2\geq L_{\alpha}$.} gives the slope of the segment $\overset\longrightarrow{AB}$ within the cone in which it lies.

The segment $\overset\longrightarrow{BF'}$ is the image of $\overset\longrightarrow{BA}$ under the map $\GG$ as shown in \cref{fig:first}; in this case it has been twisted vertically upon entering the square $S$. In \cref{fig:first}, this vertically twisted segment is shown to be cut-off at the lower edge of the square $S$. 

If this segment was not cut-off, then $\GG(\overset\longrightarrow{BA})\cap S=\overset\longrightarrow{BF'}$, and then \cref{fig:first} would be modified to include the entire segment $\overset\longrightarrow{BF'}$ in the square $S$, and we would have

\begin{align}
l_{v}(\overset\longrightarrow{BF'})=\overset\longrightarrow{BE'}=x(L_{2}+\alpha)=\beta_1 l_{v}(\gamma)(L_{2}+\alpha),
\end{align}

Let $\delta>1$ be a constant to be optimized later.

If from the above, $|\overset\longrightarrow{BE'}|=x(L_{2}+\alpha)=\beta_1 l_{v}(\gamma)(L_{2}+\alpha)>\delta l_{v}(\gamma)$, i.e.,

\begin{align}\label{eq:eq13}
   \beta_1(L_2+\alpha)>\delta,
\end{align}

we have a segment within $S$ satisfying the required property, and we would be done. Note that this $\delta$ is the same as the one used in \cref{deltareference}.

\bigskip

In case the segment $\FF(\overset\longrightarrow{BF'})$ gets cut off and $\overset\longrightarrow{BF'}\cap S=\overset\longrightarrow{BF}$, we have the following two possibilities:

\begin{enumerate}[leftmargin=*]
\item $|\overset\longrightarrow{BE}|=y(L_{2}+\alpha)>\delta\l_{v}(\gamma)$, we would be done. 

\item If not, we would have that $y\leq \frac{\delta l_{v}(\gamma)}{(L_{2}+\alpha)}$. In this case, we have further two possibilities:

\begin{enumerate}[leftmargin=*]
    \item The segment $\FF(\overset\longrightarrow{BF})$ intersects the right edge of $S$. In this case, we have, 
        \begin{align}
            l_{h}(\overset\longrightarrow{FG})>(x-y)>\beta_1 l_{v}(\gamma) -\frac{\delta l_{v}(\gamma)}{(L_{2}+\alpha)}=l_{v}(\gamma)\Big( \beta_1  -\frac{\delta}{(L_{2}+\alpha)} \Big). 
        \end{align}

        Further in that case, we have a total length of 

        \begin{align}
            l_{v}(\GG(\overset\longrightarrow{FG}))=l_{h}(\overset\longrightarrow{FG})(L_{4}+\alpha) >l_{v}(\gamma)\Big( \beta_1  -\frac{\delta}{(L_{2}+\alpha)}\Big)(L_{4}+\alpha), 
        \end{align}

    In case this was completely within the square $S$, or otherwise, we would have $\GG(\overset\longrightarrow{FG})$ intersecting the top edge of $S$. Here $L_4$ gives the slope of $\overset\longrightarrow{FG}$ within its cone.

    If we require

    \begin{align}\label{eq:eq16}
         l_{v}(\gamma)\Big( \beta_1  -\frac{\delta}{(L_{2}+\alpha)}\Big)(L_{3}+\alpha)>\delta l_{v}(\gamma),
    \end{align}

    we will be done.  Otherwise, we now have a horizontal and a vertical segment within the square. Each of these touch two adjacent sides of the square $S$, and we are done for this step.

    \item The segment $\FF(\overset\longrightarrow{BF})$ does not intersect the right edge of $S$. In this case, we ensure that $\FF(\overset\longrightarrow{BF})>\delta l_{v}(\gamma)$ and we would again be done\footnote{In this case as well, \cref{fig:first} would need to be modified.}. This follows from the following dichotomy:
    
    \begin{enumerate} 
    \item Either the segment $|\overset\longrightarrow{BE}|=l_{v}(\mathbf{F}(\overset\longrightarrow{BF}))>|\overset\longrightarrow{BD}|$ is long enough so that $l_{h}(\FF(\overset\longrightarrow{BF}))> \delta l_{v}(\gamma)$.  
    
    \item Otherwise, the segment $\overset\longrightarrow{BD}$ is sufficiently short so that consequently by virtue of the horizontal twist of $\gamma$ to the segment $\overset\longrightarrow{AC}$, the segment $|\overset\longrightarrow{JA}|=l_{h}(\overset\longrightarrow{AB})>\delta l_{v}(\gamma)$ and we would be done then itself.
    \end{enumerate}
    We make this precise here.
    
    We either have (ii) $|\overset\longrightarrow{JA}|=l_{v}(\overset\longrightarrow{AB})>\delta l_{v}(\gamma)$ in which case we are done. Otherwise, with some basic geometry, we verify that we have (i):

    \begin{align}
        |\overset\longrightarrow{BE}|\geq |\overset\longrightarrow{BD}|=l_{v}(\gamma) - \frac{|\overset\longrightarrow{JA}|}{L_{1}+\alpha}\geq l_{v}(\gamma)\Big( 1- \frac{\delta}{L_{1}+\alpha} \Big).
    \end{align}
    
    In this case, it is enough to ensure

    \begin{align}\label{eq:eq18}
    (L_{3}+\alpha)l_{v}(\gamma)\Big( 1- \frac{\delta}{L_{1}+\alpha} \Big)>\delta l_{v}(\gamma),
    \end{align}

    since we have $l_{h}(\FF(\overset\longrightarrow{BF}))=(L_{3}+\alpha)l_{v}(\gamma)\Big( 1- \frac{\delta}{L_{1}+\alpha} \Big)$, which forces \\  $l_{h}(\FF(\overset\longrightarrow{BF}))>\delta l_{v}(\gamma)$.
    
\end{enumerate}

Here, $L_{3}$ is the slope of the segment $\overset\longrightarrow{BF}$ within the cone in which it lies.

\bigskip

\end{enumerate}

Further, assume that we haven't already found a segment above, that has vertical or horizontal length at least $l_v(\gamma)$ (otherwise, we would be done). 

In this case, we ensure that either $\FF(\li{GH})> \delta l_{v}(\gamma)$ or that  $\FF(\li{GH})$ touches the left end of the square $S$. Note that $l_h(\li{FG})\leq \delta l_{v}(\gamma)$, and thus that $|L_4| l_h(\li{FG})\leq |L_4| \delta l_{v}(\gamma)$. Further we clearly have that $|\li{KP}|>l_{v}(\gamma)$, and thus $l_{v}(\li{GH})>l_{v}(\gamma)-|L_4| \delta l_{v}(\gamma)=(1-|L_4|)l_v (\gamma)$. In this case, either we have, 

\begin{align*}
\FF(\li{GH})>(L_{5}+\alpha)(1-\delta |L_4|)l_{v}(\gamma),
\end{align*}

or that $\FF(\li{GH})$ touches the left edge of $S$. With the above, we note that the sufficiency of $(L_{5}+\alpha)(1-\delta |L_4|)l_{v}(\gamma)>\delta l_v (\gamma)$ is equivalent to \cref{eq:eq18} for our purpose of optimization of $\alpha$ later on, and this is not a new restriction on $\alpha$. Further on, if $\l_{v}(\li{HG})<\delta l_{v}(\gamma)$, then $|\li{HK}|<\delta |L_5| l_v (\gamma)$ where $L_5$ is the slope of the segment $\li{HG}$ within its cone. Further, we have assumed that the horizontal length of $S$ which is $|\li{EP}|$, satisfies $|\li{EP}|> |\li{JA}|=\beta_1 l_{v}(\gamma)$. In this case, we have that $l_{h}(\li{MH})>\beta_1 l_{v}(\gamma)-\delta |L_5| l_v (\gamma) $, and then 

\begin{align*}
\GG(\li{HM})>(L_6 +\alpha)(\beta_1 -\delta |L_5|)l_v (\gamma),
\end{align*}

or that $\GG(\li{HM})$ touches the bottom edge of the square $S$. Here $L_6$ is the slope of the segment $\li{HM}$ within its cone. Again with the above, we will note that the sufficiency of the condition $(L_6 +\alpha)(\beta_1 -\delta |L_5|)l_v (\gamma)>\delta l_v(\gamma)$ is equivalent to \cref{eq:eq16} above.

The same dynamics repeats further on, till we get either a vertical or a horizontal segment that has length greater than $\delta l_v (\gamma)$ and we are done, or we get segments with successively the slopes 

\begin{align}
L_{i+1}=-\frac{1}{L_i +\alpha},
\end{align}

and we note that as $i\to \infty$, these slopes $L_i$ converge to $L_\alpha$.

One way to see this is to consider the function $f(x)=x(x+\alpha)$ when restricted to values of $L_\alpha \leq x\leq 0$. This function is seen to be increasing whenever $x\geq -\alpha/2$ which is true in the above domain for $x$ and that $f(0)=0, f(L_\alpha)=-1$ . In this case, if for a given $x\in [L_\alpha,0]$ we have $f(x)=x(x+\alpha)=-t_1$ for $0\leq t\leq1$, then $-(1/(x+\alpha))=x/t<x $, and plugging in $1/(x+\alpha)$ in place of $x$ and doing the analysis again yields $t_2>t_1$ with which the above repeats again, with the $t_i \to 1$ as $i \to \infty$, and that finally the sequence of values converge to $L_\alpha$.

Thus in the limit, we would have sequences of four connected segments that converge as shown in part (a) of \cref{fig:figsec}. We note that given any such $S$, there is a unique way to have four connected segments in the above manner so that each has slope $L_\alpha$ within its cone, and our segments converge to this limiting sequence.

\bigskip

\subsection{Part 2}
Next we consider the cases where a segment $\overset{\sim}{\gamma}$ of horizontal length at least $d$, enters inside the square $S$. There are two possibilities: 

\begin{enumerate}
    \item The segment $\overset{\sim}{\gamma}$ touches the left edge of $S$ and lies above the line on which the rational orbit lies . In this case, we ensure that the vertical length of the segment $I_{1}\cup I_{2}'$ is long enough so that there is a sufficiently long vertical distance below the line of the rational orbit to the bottom layer of $S$. After this we employ an argument similar to the one employed earlier for the case of the segment $I_{4}$, as will be clear from the argument that follows.

    \item  This segment $\overset{\sim}{\gamma}$ touches the right edge of $S$. In this case this segment could be arbitrarily close to the top or bottom edge of S.
\end{enumerate}

\bigskip

Both of these above situations involve arguments analogous to the ones presented in the previous subsection. In addition, when in the second case the segment under consideration is arbitrarily close to the top edge of $S$, we must deal with the contingency that within $S$ we do not end up with either a long enough segment under two successive iterations of $\FF,\GG$, nor a combined segment consisting of both a v-segment and an h-segment. We then iterate further under the maps $\FF,\GG$ successively outside the square $S$ to eventually find a long enough segment. 

\bigskip

Consider \cref{fig:second}(b), where $F$ now acts vertically and $G$ acts horizontally, and where $\l_{v}(\overset\longrightarrow{AB})=\beta_{2} l_{v}(\gamma)$ with $0< \beta_{2}<1$ to be later optimized to be as small as possible. Let $L_{5}$ be the slope of the segment $\overset\longrightarrow{AB}$ within it's cone. Under the map $\Phi$, we have either of the following three cases: 

\begin{enumerate}[leftmargin=*]
    \item $l_{h}(\GG(\overset\longrightarrow{AB}))>\delta l_{v}(\gamma)$ and $\GG(\overset\longrightarrow{AB})$ does not touch the right edge of $S$, in which case we are done.

    \item $\GG(\overset\longrightarrow{AB})$ touches the right edge of $S$ but we still have $l_{h}(\GG(\overset\longrightarrow{AB})\cap S)=l_{h}(\overset\longrightarrow{BC})>\delta l_{v}(\gamma)$, in which case we are also done. With some elementary geometry, referring to \cref{fig:second}(b), we have that 
    
    \begin{align}
    \frac{|\overset\longrightarrow{CX}|}{l_{h}(\overset\longrightarrow{BC})}=\frac{1}{(L_{5}+\alpha)}.
    \end{align}

    Thus in this instance, we must have $l_{v}((\overset\longrightarrow{BC}))= |\overset\longrightarrow{CX}|> \frac{\delta l_{v}(\gamma)}{(L_{5}+\alpha)}$. 

    \item $\GG(\overset\longrightarrow{AB})$ touches the right edge of $S$, but now 
    
    \begin{align}\label{eq:eqimp}
    l_{v}((\overset\longrightarrow{BC}))\leq \frac{\delta l_{v}(\gamma)}{(L_{5}+\alpha)}
    \end{align}
    
    In this case, we must have either of the two following cases:
    
    \begin{enumerate}[leftmargin=*]
        \item 
    $\FF(\overset\longrightarrow{BC})$ touches the top edge in \cref{fig:second}(b) \footnote{Note that, as drawn, the left edge of $S$ in \cref{fig:second} is the top edge of $S$ in \cref{fig:second}, and so forth. The horizontal lengths in \cref{fig:second}(a) also become the vertical lengths in \cref{fig:second}(b), and so forth. The references to the figure and edges, and the horizontal and vertical directions will be made clear from context.}. In this case, using \cref{eq:eqimp}, we have clearly:
    \begin{align}
        |\overset\longrightarrow{CY})|=l_{v}(\FF(\overset\longrightarrow{BC})))\geq l_{v}(\overset\longrightarrow{AB}))-l_{v}(\overset\longrightarrow{BC}))>\beta_2 l_{v}(\gamma)-\frac{\delta l_{v}(\gamma)}{(L_{5}+\alpha)}=l_{v}(\gamma)\Big( \beta_2 -\frac{\delta}{(L_{5}+\alpha)} \Big)
    \end{align}

    In this case, if $l_{v}(\FF(\li{BC}))> l_{v}(\gamma)$, we are done. Let $L_{7}$ be the slope of the segment $\overset\longrightarrow{CE}$ within it's cone. In this case, we then have that either,
    
    \begin{enumerate}[leftmargin=*]
    \item $    \GG(\overset\longrightarrow{CE}))$ touches the left edge of $S$ in \cref{fig:second}(b). In this case, we further have to ensure that $\FF(\li{EJ})$ either touches the bottom edge of $S$ in \cref{fig:second}(b), or else that we have a segment $\FF(\li{EJ})$ whose vertical length is greater than $\delta l_v (\gamma)$. Call the slope of the segment $\li{EJ}$, $L_7$ within its cone. Also it is clear that in \cref{fig:second}(b), analogous to the previous case, the horizontal length of $S$ is at least $l_{v}(\gamma)$, since the original segment $\overset\longrightarrow{A'C'}$ lies inside the square $S_1$ as in \cref{fig:second}(a). Also since $l_{v}(\li{CE})<\delta l_{v}(\gamma)$, we have $l_{h}(\li{CE})<|L_8|\delta l_{v}(\gamma)$ where $L_8$ is the slope of the segment $\li{CE}$ within its cone. Thus we have $l_{h}(\li{EJ})> l_{v}(\gamma)(1-|L_8|\delta)$, and thus further, that $l_{v}(\FF(\li{EJ}))>(\alpha+ L_9)  l_{v}(\gamma)(1-|L_8|\delta)$, where $L_9 =1 /(L_8 +\alpha)$ is the slope of the segment $\li{JK}$ within its cone.

    Thus, the following condition

    \begin{align}\label{eq:eq23}
        (\alpha +L_9)(1-|L_8|\delta)>\delta,
    \end{align}

    is sufficient to get a long enough segment $\li{JK}$ or that $\li{JK}$ touches the bottom edge of $S$. In case this segment is cut off by the bottom edge of $S$ and is of length less than $\delta l_{v}(\gamma)$. Further on, note that the vertical length of $S$ is at least $\beta_2 l_{v} (\gamma)$ and that $l_{v}(\li{JK})>\beta_{2}l_{v}(\gamma)-|L_8|\delta l_{v}(\gamma)=(\beta_{2} -|L_8|\delta)l_{v}(\gamma)$. Thus we have that either $\GG(\li{JK})$ is cut off by the right edge of $S$ or that the horizontal length $l_{h}(\GG{(\li{JK}})$ is at least $(L_{9} +\alpha)(\beta_2 -|L_8|\delta)l_{v}(\gamma)$. Thus the following condition is sufficient,

    \begin{align}\label{eq:eq24}
        (L_9 +\alpha)(\beta_2 -|L_8|\delta)>\delta,
    \end{align}

    to get a long enough segment if $\GG{(\li{JK})}$ doesn't get cut off by the right edge of $S$. 

    Further on, with requirements of the same form as the above two conditions \cref{eq:eq23,eq:eq24}, we can ensure that successively  we have segments such as in Figure 2(a) or that we have a segment of horizontal or vertical length at least $l_{v}(\gamma)$ in which case we are done. 
 
    \item Otherwise, $\GG(\overset\longrightarrow{CE})$ does not intersect the left edge of $S$, in which case it will be enough to ensure that:

    \begin{align}\label{eq:eq22}
    l_{h}(\GG(\overset\longrightarrow{CE})))\geq(L_{7}+\alpha) l_{v}(\gamma)\Big( \beta_2 -\frac{\delta}{(L_{5}+\alpha)} \Big)  >\delta l_{v}(\gamma).
    \end{align}

    We note that the requirement from this above equation is identical to the requirement from \cref{eq:eq24}.

    \end{enumerate}

    \item $\FF(\overset\longrightarrow{BC})$ lies entirely within $S$, but does not touch the top edge of $S$ in \cref{fig:second}(b), in which case we have the following two possibilities:

    \begin{enumerate}[leftmargin=*]

     \item We have $l_{v}(\FF(\overset\longrightarrow{BC}))> \delta l_{v}(\gamma)$ in which case we are done.

    If the slope of the segment $\overset\longrightarrow{BC}$ is $L_{6}$ within its cone \footnote{In this case this cone is oriented horizontally.}, then referring to \cref{fig:second}(b), it happens that,  

\begin{align}
    l_{v}(\FF(\overset\longrightarrow{BC})))=(L_{6}+\alpha)l_{h}(\overset\longrightarrow{BC})=(L_{6}+\alpha)|\overset\longrightarrow{BX}|>\delta l_{v}(\gamma)
\end{align}

     \item On the other hand if we have  $l_{v}(\FF(\overset\longrightarrow{BC}))\leq \delta l_{v}(\gamma)$, then:

\begin{align}
    |\overset\longrightarrow{BX})|=l_{h}(\overset\longrightarrow{BC})\leq \frac{\delta l_{v}(\gamma)}{(L_{6}+\alpha)}
\end{align}

Further, in this case, we also have:

\begin{align}\label{CX}
    |\overset\longrightarrow{CX}|=|L_{6}|\cdot l_{h}(\overset\longrightarrow{BC})\leq \frac{|L_{6}|\delta l_{v}(\gamma)}{(L_{6}+\alpha)}
\end{align}

    Thus in this case, we have 

    \begin{align}
        \eta=\frac{l_{v}(\overset\longrightarrow{AB})}{|\overset\longrightarrow{CX}|}= \frac{\beta_2 l_{v}(\gamma)}{|\overset\longrightarrow{CX}|}\geq \frac{\beta_2 (L_{6} +\alpha)}{|L_{6}|\delta}
    \end{align}

     In this case, we consider the orbit $G^{k}(C)|_{k=1}^{\infty}$ of the point $C$ under the horizontal twist $G$, as shown in the lift of the track in \cref{fig:second}(c).

     \bigskip
    Note that as a result of one horizontal twist $G$ in \cref{fig:second}(c), we have $C=\GG(C')$ and in particular, $d(C',C)<D_{1}$. Under a further horizontal twist, the point $C$ moves again a distance $d(C',C)$ which is also at most $D_{1}$. This distance is also equal to $|\overset{\longrightarrow}{CC'}|=d(C,\GG(C))=(\beta_2 l_{v}(\gamma)\alpha)/\eta$ . Also, since $D_{2}> 0$, it is easy to see that $\GG(C)$ cannot lie to the right of $RE_{2}$.

    Also, using Equations (27) and (28), we see that

    \begin{align}
        |\overset\longrightarrow{CC'}|=|\overset\longrightarrow{BX}|+|L_{5}||\overset\longrightarrow{CX}|\leq \frac{\delta l_{v}(\gamma)}{(L_{6}+\alpha)}\Big( 1+|L_{5}L_{6}| \Big)
    \end{align}

    Further, we have that 

    \begin{align}
        |\overset\longrightarrow{AA''}|=|\overset\longrightarrow{BX}|+|L_{5}|\beta_2 l_{v}(\gamma)\leq \frac{\delta l_{v}(\gamma)}{L_{6}+\alpha}+|L_{5}|\beta_2 l_{v}(\gamma)
    \end{align}

   \begin{enumerate}
       \item 
   Suppose that $\GG(C)$ lies in $S_2$, in between $LE_{2}$ and $RE_{2}$ (including possibly on either of the edges). Thus, in particular, $d(RE_{1},LE_{2})<|CC'|$. Also note that in this case, $d(A,\GG(A))=\alpha\beta_2 l_v (\gamma)$.
    
    In this case the Figure 6(c) is not to scale and we would have $A'$ lying to the right of $LE_2$, and unless we get a horizontal segment through $S_2$ and we are done, we have using Equations (27) and (28):
    \begin{align}
     l_{h}(\GG(\overset\longrightarrow{C'A})\cap S_{{2}})=\alpha \beta_2 l_{v}(\gamma)-d(RE_{1},LE_{2})-|\overset\longrightarrow{AA''}| \\ \implies l_{h}(\GG(\overset\longrightarrow{C'A})\cap S_2) > \alpha \beta_2 l_{v}(\gamma)-|\overset\longrightarrow{CC'}|-|\overset\longrightarrow{AA''}|\\ \implies l_{h}(\GG(\overset\longrightarrow{C'A})\cap S_2) > l_{v}(\gamma)\Big( \alpha\beta_2 -\frac{\delta(2+|L_{5}L_{6}|)}{L_{6}+\alpha}-\beta_2 |L_{5}| \Big)
    \end{align}

    \begin{figure}
\centering
    \begin{tikzpicture}
        \draw[very thick] (-3,-1.2) to (-3,3);
        \draw[very thick] (-1.2,-1.2) to (-1.2,3);
        \draw[very thick] (2.7,-1.2) to (2.7,3);
        \draw[very thick] (4.5,-1.2) to (4.5,3);
        \draw[very thick] (-3.4,2.8) to (4.9,2.8);
        \draw[very thick] (-3.2,-0.5) to (4.7,-0.5);
        \draw[dashed, <->] (-3.4,-0.5) to (-3.4,2.8);
        \draw[dashed, <->] (-4.4,0.5) to (-4.4,2.7);
        \draw[dashed, <->] (3.2,-1.1) to (4.5,-1.1);
        \node [above] at (-3.7,0.6) {$D_{1}$};
        \draw[very thick] (-2,0.5) to (-2.2,2.7);
        \node [below] at (-2,0.2) {$S_1$};
        \node [below] at (4,-1.1) {$l_{h}(\textbf{F}(\gamma)\cap S)$};
        \node [below] at (3.45,0.2) {$S_2$};
        \node [left] at (-5.5,1.7) {$(a)$};
        \node [left] at (-4.3,1.9) {$l_{v}(\gamma)$};
        \draw[very thick] (3.2,2.4) to (4.5,2.6);
        \draw (4.38,2.8) to (4.5,2.6);
        \draw (4.38,2.8) to (2.7,2.65);
        \draw (3.9,-0.5) to (2.7,2.65);
        \node [below] at (3.2,2.4) {$A$};
        \node [right] at (4.45,2.6) {$B$};
        \node [above] at (4.35,2.8) {$C$};
        \node [left] at (2.76,2.65) {$E$};
        \node [below] at (3.9,-0.5) {$J$};
        \draw (3.9,-0.5) to (4.5,0.2);
        \node [right] at (4.5,0.2) {$K$};
        \draw[->] (6.08,1) to (6.08,0.5);
        \draw[->] (6.08,1) to (6.58,1);
        \node [left] at (6.08,0.6) {$\GG$};
        \node [above] at (6.58,1) {$\FF$};

        \draw[very thick] (-2,-3) to (3,-3);
        \draw[very thick] (-2,-5.8) to (3,-5.8);
        \draw[very thick] (-1.5,-2.8) to (-1.5,-6.0);
        \draw[very thick] (2.3,-2.8) to (2.3,-6.0);
        \node [left] at (-5.5,-4.7) {$(b)$};
        \node [left] at (0.5,-4.7) {$S_1$};
        \node [above] at (0.5,-2.6) {$D_1$};
        \draw[dashed, <->] (-1.5,-2.5) to (2.3,-2.5);
        \draw[dashed] (4.9,-3.5) to (1.7,-3.5);
        \draw[dashed] (3.6,-5.8) to (1.8,-5.8);
        \draw[very thick] (1.8,-3.5) to (2.0,-5.8);
        \draw (2.0,-5.8) to (5,-3.5);
        \draw (2.08,-3) to (2.3,-5.65);
        \draw (2.08,-3) to (-1.5,-4.3);
        \draw (-1.0,-5.8) to (-1.5,-4.3);
        \node [below] at (-1.0,-5.8) {$K$};
        \node [left] at (1.9,-3.7) {$A$};
        \node [below] at (1.85,-5.75) {$B$};
        \node [right] at (2.3,-5.45) {$C$};
        \node [left] at (2.1,-5.45) {$C'$};
        \draw[dashed] (2.3,-5.57) to (2.0,-5.57);
        \node [right] at (2.20,-5.99) {$X$};
        \node [above] at (2.02,-3) {$E$};
        \node [below] at (2.45,-2.9) {$Y$};
        \node [left] at (-1.5,-4.3) {$J$};
        \draw[dashed, <->] (3.6,-3.5) to (3.6,-5.8);
        \node [right] at (4.8,-3.7) {$A'$};
        \node [right] at (2.2,-3.7) {$A''$};
        \node [right] at (3.6,-4.9) {$\beta_2 l_{v}(\gamma)$};
        \draw[->] (-3.08,-4) to (-3.08,-3.5);
        \draw[->] (-3.08,-4) to (-2.58,-4);
        \node [left] at (-3.08,-3.7) {$\FF$};
        \node [below] at (-2.58,-4) {$\GG$};

        \draw[very thick] (-3,-8) to (11.5,-8);
        \draw[very thick] (-3,-9.5) to (11.5,-9.5);
        \node [left] at (-5.5,-9) {$(c)$};
        \draw[very thick] (-2.5,-7.8) to (-2.5,-9.7);
        \draw[very thick] (0,-7.8) to (0,-9.7);
        \draw[very thick] (3,-7.8) to (3,-9.7);
        \draw[very thick] (5.5,-7.8) to (5.5,-9.7);
        \draw[very thick] (8.5,-7.8) to (8.5,-9.7);
        \draw[very thick] (11,-7.8) to (11,-9.7);
        \node [above] at (-1.2,-8.5) {$S_1$};
        \node [above] at (4.25,-8.5) {$S_2$};
        \node [above] at (9.75,-8.5) {$S_3$};
        \draw[very thick] (-0.2,-9.5) to (-0.3,-8.3);
        \draw (2,-8.35) to (-0.2,-9.5);
        \node [right] at (2,-8.35) {\footnotesize $A'$};
        \node [left] at (-0.2,-8.35) {\footnotesize $A$};
        \node [right] at (0,-9.4) {\footnotesize $C$};
        \node [left] at (-0.15,-9.4) {\footnotesize $C'$};
        \node [below] at (-0.2,-9.42) {\footnotesize $B$};
        \draw[dashed] (-0.2,-9.42) to (0,-9.42);
        \draw [fill] (0,-9.42) circle [radius=0.05];
        \draw [fill] (1.6,-9.42) circle [radius=0.05];
        \draw [fill] (3.2,-9.42) circle [radius=0.05];
        \draw [fill] (4.8,-9.42) circle [radius=0.05];
        \draw [fill] (6.4,-9.42) circle [radius=0.05];
        \draw [fill] (8.0,-9.42) circle [radius=0.05];
        \draw [fill] (9.6,-9.42) circle [radius=0.05];
        \draw [fill] (11.2,-9.42) circle [radius=0.05];
        \draw[dashed, <->] (-2.5,-9.8) to (0,-9.8);
        \draw[dashed, <->] (0,-9.8) to (3,-9.8);
        \draw[dashed, <->] (3,-9.8) to (5.5,-9.8);
        \draw[dashed, <->] (5.5,-9.8) to (8.5,-9.8);
        \draw[dashed, <->] (8.5,-9.8) to (11,-9.8);
        \node [below] at (-1.25,-9.76) {$D_1$};
        \node [below] at (1.75,-9.76) {$D_2$};
        \node [below] at (4.25,-9.76) {$D_1$};
        \node [below] at (7.25,-9.76) {$D_2$};
        \node [below] at (10.0,-9.76) {$D_1$};
        \node [above] at (0,-8) {$RE_1$};
        \node [above] at (5.4,-8) {$RE_2$};
        \node [above] at (3.2,-8) {$LE_2$};
        \node [above] at (8.4,-8) {$LE_3$};
        \node [above] at (10.9,-8) {$RE_3$};


    \end{tikzpicture} 
    \caption{In part (a), we lift to $\R$ and denote the successive lifts of the square $S$ as $S_1$ and $S_2$ as in the figure. Part (b) shows an enlarged picture of the square $S_2$, rotated by $\pi/2$, and the case where under iterations of the twist, the segment $\overset\longrightarrow{BC}$ has a long enough twist in the square $S$. In part (c) the black dots represent the points of the orbit of $C$ under the successive `horizontal' twists in this figure, with the three successive squares shown as $S_1,S_2,S_3$ in this figure. (The part (c) is not up to scale.)}
    \label{fig:second}
\end{figure}
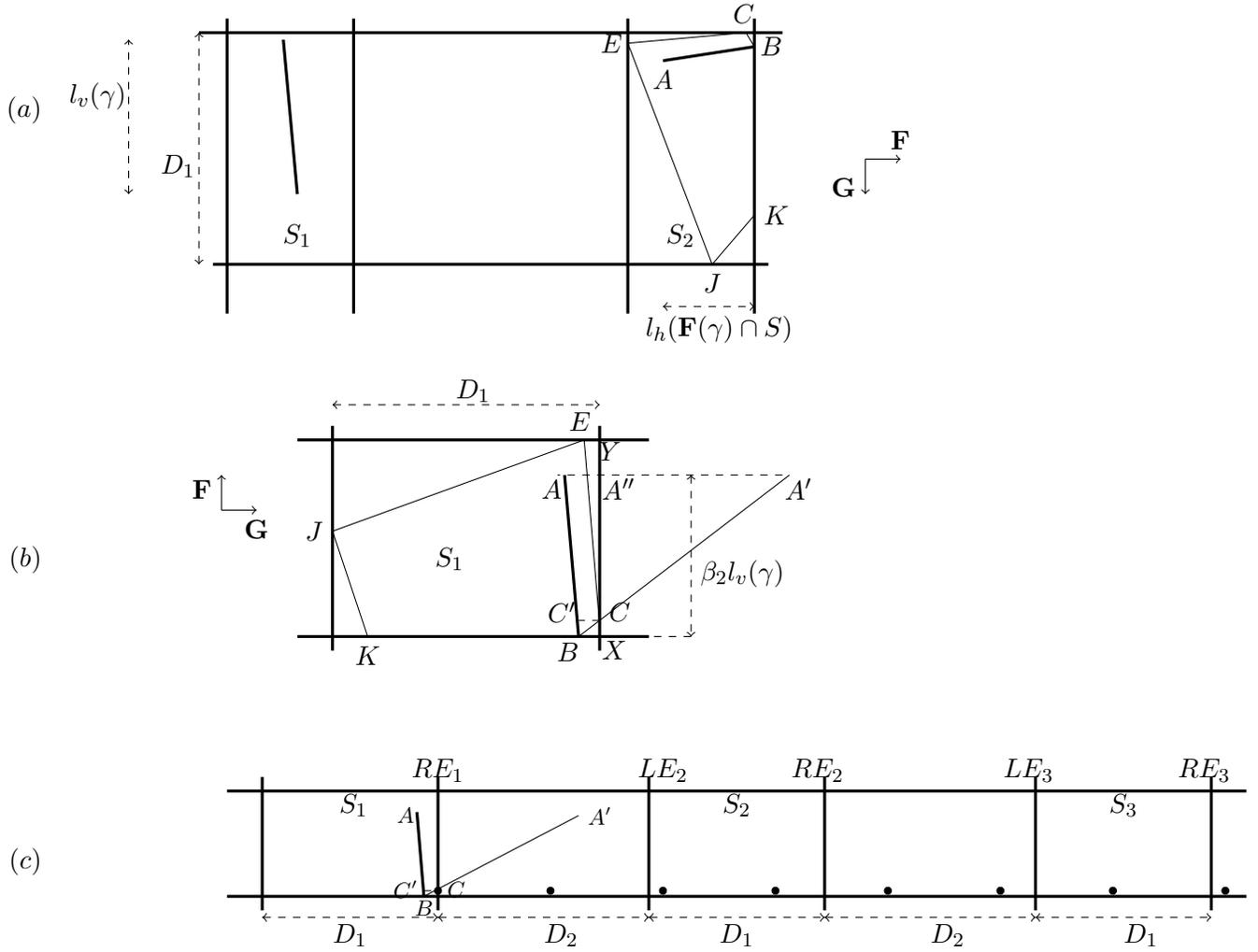

    Thus we will be done if:

    \begin{align}\label{eq:eq32}
        \Big( \alpha\beta_2 -\frac{\delta(2+|L_{5}L_{6}|)}{L_{6}+\alpha}-\beta_2 |L_{5}| \Big)>\delta 
    \end{align}

    \bigskip

   \item Now suppose that $\GG(C)$ lies between $RE_{1}$ and $LE_{2}$.
    
    In this case, suppose there are $m\geq 1$ many points of the orbit of $G^{k}(C)$ in the region between $RE_{1}$ and $LE_{2}$(not including these two edges). Let $t$ be the least integer with $1\leq t\leq m$, for which $G^{t}(A)\in S_2$. For this we look at \cref{fig:third}:

    In this case, if $|\overset\longrightarrow{Q_1Q_2}|> \delta l_{v}(\gamma)$ we are done, otherwise we assume that $|\overset\longrightarrow{Q_1Q_2}|\leq \delta l_{v}(\gamma)$, and we have with some elementary geometry that 

    \begin{align}
        \frac{r_1}{r_2}=\frac{|\li{Q_1 Q_2}|}{|\li{P_1 P_2}|}\leq \frac{\delta l_{v}(\gamma)}{t(1- \frac{1}{\eta})\beta_2\alpha l_{v}(\gamma)-|L_5|(1- \frac{1}{\eta})\beta_2 l_{v}(\gamma)-\delta l_{v}(\gamma)} \\ \leq \frac{\delta l_{v}(\gamma)}{(1- \frac{1}{\eta})(\alpha-|L_5|)\beta_2 l_{v}(\gamma)-\delta l_{v}(\gamma)}
    \end{align}

   Note that $|L_5|(1- \frac{1}{\eta})\beta_2 l_{v}(\gamma)$ is the length of the horizontal projection of $\li{AC'}$. We also have, looking at \cref{fig:third}, that $r_1 +r_2= (1- 1/\eta)\beta_2 l_{v}(\gamma)$. With an elementary calculation, we can verify that:

    \begin{align}
        r_1\leq (r_1 +r_2) \frac{\delta l_{v}(\gamma)}{(1- \frac{1}{\eta})(\alpha -|L_5|)\beta_2 l_v (\gamma)}= \frac{\delta l_{v}(\gamma)}{\alpha-|L_5|}
    \end{align}

    Thus, we have 

    \begin{align}\label{eq:eq39}
        |\overset\longrightarrow{MN}|= \beta_2 l_{v}(\gamma)-r_1\geq l_{v}(\gamma)\Big( \beta_2 -\frac{\delta}{\alpha-|L_5|} \Big)
    \end{align}
When we have $t=1$, we note that by hypothesis there exists at least one more point $G^{2}(C')=G(C)$ that lies in between $RE_1$ and $LE_2$.

In this case, in the 2nd iteration, we clearly have either a segment of length $\alpha |\overset\longrightarrow{MN}|=\alpha l_{v}(\gamma)\big( \beta_2 -\delta/(\alpha-|L_5|)\big) $ within the square $S_2$ or there is a horizontal segment within $S_2$. Thus we would be done if we require:

    \begin{align}\label{eq:eq36}
        \alpha \Big( \beta_2 -\frac{\delta}{(\alpha-|L_5|)}\Big) > \delta.
    \end{align}

\bigskip
Now if we have $t\geq 2$, which was defined above as the least integer so that $G^{t}(A)\in S_2$, then we need to ensure that there exists at least one further point $G^{t+1}(C')$ that lies in between $RE_1$ and $LE_2$(including possibly on the edge $LE_2$). Then, by the argument above, in the $(t+1)$'th iteration, this would force again a segment of length at least $\delta l_{v}(\gamma)$ inside $S_2$ or a horizontal segment through $S_2$.

In order to achieve this, we first note that it is enough to consider the case of $t=2$. This is because for $t\geq 1$, the horizontal projection of the segment joining $G^t(C')$ and $G^t(A)$ is an increasing function of $t$.

 For the case of $t=2$, we first note that $G(A)$ lies between $RE_1$ and $LE_2$. In this case, it will be enough to ensure that $G^{3}(C')$ lies to the left of $G(A)=A'$. In that case, one can easily see that in the $t=3$ iteration we would ensure the requisite segment inside $S^2$. 

To ensure this, it is enough to compare $(\alpha-|L_5|)|\li{CX}|+2\alpha|\li{CX}|=(3\alpha -|L_5|)|CX|$ and, $(\alpha -|L_5|)|\beta_2 l_{v}(\gamma)$ . Further, using the upper bound from \cref{CX} on $|\li{CX}|$, it is enough to require that,
\begin{align}\label{new1}
     (\alpha - |L_5|)\beta_2 l_{v}(\gamma)\geq (3\alpha -|L_5|)|\li{CX}|.
\end{align}
For this, it is enough to require that,
\begin{align}\label{new2}
    (\alpha - |L_\alpha|)\beta_2 \geq 3\alpha \frac{|L_\alpha|\delta }{(\alpha-|L_\alpha|)} \Leftrightarrow \beta_2 \geq \frac{3\alpha |L_\alpha|\delta}{(\alpha- |L_\alpha|)^{2}}
\end{align}

    \end{enumerate}
\end{enumerate}
    \end{enumerate}
\bigskip
    
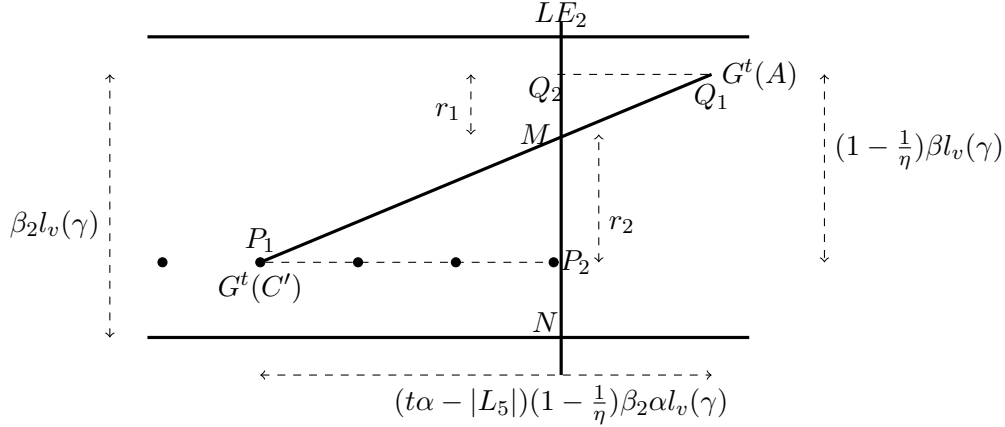
\begin{figure}
\centering
    \begin{tikzpicture}
        \draw[very thick] (1.5,2.2) to (1.5,-2.5);
        \draw[very thick] (-4,2) to (4,2);
        \draw[very thick] (-2.5,-1) to (3.5,1.5);
        \draw[very thick] (-4,-2) to (4,-2);
        \draw[dashed] (1.45,1.5) to (3.5,1.5);
        \draw[dashed] (1.47,-1) to (-2.5,-1);
        \draw[dashed,<->] (5,1.5) to (5,-1);
        \draw [fill] (-2.5,-1) circle [radius=0.06];
        \draw [fill] (-3.8,-1) circle [radius=0.06];
        \draw [fill] (-1.2,-1) circle [radius=0.06];
        \draw [fill] (0.1,-1) circle [radius=0.06];
        \draw [fill] (1.4,-1) circle [radius=0.06];
        \node [below] at (-2.5,-1) {$G^{t}(C')$};
        \node [right] at (3.5,1.5) {$G^{t}(A)$};
        \node [above] at (1.5,2.0) {$LE_2$};
        \node [right] at (5,0.5) {$(1- \frac{1}{\eta})\beta_2 l_{v}(\gamma)$};
        \draw[dashed, <->] (-2.5,-2.5) to (3.5,-2.5);
        \node [below] at (1.5,-2.5) {$(t\alpha -|L_5|)(1- \frac{1}{\eta})\beta_2 l_{v}(\gamma)$};
        \draw[dashed, <->] (2,-1) to (2,0.7);
        \node [right] at (2,-0.5) {$r_2$};
        \node [right] at (1.35,-1) {$P_2$};
        \node [above] at (-2.5,-1) {$P_1$};
        \node [left] at (1.67,1.3) {$Q_2$};
        \node [below] at (3.5,1.5) {$Q_1$};
        \draw[dashed, <->] (-4.5,-2) to (-4.5,1.5);
        \node [left] at (-4.5,-0.5) {$\beta_2 l_{v}(\gamma)$};
        \node [left] at (1.5,0.7) {$M$};
        \node [left] at (1.6,-1.8) {$N$};
        \node [left] at (0.3,1.0) {$r_1$};
        \draw[dashed, <->] (0.3,1.5) to (0.3,0.7);

    \end{tikzpicture} 
    \caption{The case outlined when $G^{t}(\overset\longrightarrow{AC})$ intersects $S_2$ for the first time, when considering the dynamics outlined in \cref{fig:second}(c). The point $P_1$ coincides with $G^{t}(C')$ while the point $Q_{1}$ coincides with $G^{t}(A)$.}
    \label{fig:third}
\end{figure}
 \end{enumerate}
We also have remaining the case where a segment of length at least $d$ enters the square $S$ and touches the left edge, as opposed to the above case where the segment $d$ touched the right edge of $S$. But here the situation is basically identical to the one already considered, and the above analysis covers this case, and the same parameter $\beta_{2}$ used above suffices for this case. Thus the same uniform lower bound of $\beta_2 l_{v}(\gamma)$ can be used for the two lower bounds on $d$ coming from these two cases in the subsequent analysis. 

\subsection{Part 3} We have the following four contingencies outlined in the beginning of Section 2. 

\subsubsection{Case a.} 

\begin{enumerate}
    \item A horizontal segment through $S$ belonging to $\FF\circ\Phi^{m_0}(\gamma^{u}(x))$. In this case we get either a vertical segment through $S$ belonging to $\Phi^{m_0 +1}(\gamma^{u}(x))$ or otherwise if there is also a horizontal segment through $S$ belonging to $\Phi^{-n_0}(\gamma^{s}(y))$ then we can guarantee a point of intersection between segments belonging to $\Phi^{m_0 +1}(\gamma^{u}(x))$ and $\Phi^{-n_0-1}(\gamma^{s}(y))$ in which case we are done.


\begin{enumerate}[leftmargin=*]
    \item[$\mathcal{z}_1:$] All the variables used are defined in the diagram in  \cref{fig:inksone3}(a). We have a modification of the argument from the previous section, and a situation analogous to part (c) of \cref{fig:second} is dealt with here. Initially with upper bounds on $l_1,l_2$, we set up analogues of \cref{eq:eq32} and \cref{eq:eq36}, in order to ensure a vertical segment through $S$. When either one of $l_1$ or $l_2$ violates this upper bound, we then follow a modified argument.

    Consider a parameter $\kappa$ to be determined later, so that both $l_1,l_3\leq \kappa l$. We also note from the geometry that $l_2\leq |L_{\alpha}|l$.

    Using elementary geometry, we will consider the segment $\li{BA}$ and it's iteration under $\mathbf{G}$. We note that $|CC'|\leq l_1(1+ |L_\alpha|^{2})$, and $|\li{AA''}|=l_1 +l_2\leq (\kappa +|L_\alpha|)l$. Note that the analog of \cref{eq:eq32}, when $\GG(C)$ again lies within a lift of $S$ itself\footnote{In reference to \cref{fig:second}, this means that $\GG(C)\in S_{2}$.}, gives us

    \begin{equation*}
        l_{h}(\GG(\li{C'A})\cap S)\geq \alpha l - |\li{CC'}|-|\li{AA''}|\geq \alpha l -l_1 (1+|L_{\alpha}|^{2})-(\kappa +|L_\alpha|)l,
    \end{equation*}

    and thus,

    \begin{equation}
        l_{h}(\GG(\li{C'A})\cap S)\geq \alpha l -\kappa (1+|L_\alpha|^{2})l -(\kappa +|L_\alpha|)l=l(\alpha -\kappa(1+|L_\alpha|^{2})-(\kappa +|L_{\alpha}|))
    \end{equation}

In this case, if we require that 

\begin{equation}\label{eq:eq42}
    (\alpha -\kappa(1+|L_\alpha|^{2})-(\kappa +|L_{\alpha}|))>(2\kappa +|L_{\alpha}|)\Rightarrow \alpha>4\kappa + \kappa|L_\alpha|^{2} +2|L_\alpha|,
\end{equation}

then from the above, we get that 

\begin{equation*}
     l_{h}(\GG(\li{C'A})\cap S)\geq l_1 +l_2 +l_3.
\end{equation*}

Thus in fact we have a vertical segment through $S$. Now the argument goes onto the item 2 below.

The second case (B) concerns the analogs of \cref{eq:eq39,eq:eq36,new1,new2}. First, noting that $L_{\li{AB}}$ is the slope of the segment $\li{AB}$ in the cone in which it lies, we basically repeat the argument prior to \cref{eq:eq39,eq:eq36}. If we have a segment of length $l_1 +l_2 +l_3$ already inside $S_2$ then we are done. Otherwise using the same notation as in \cref{fig:third}, we would have, 
\begin{align}
    \frac{r_1}{r_2}\leq \frac{l_1 +l_2 +l_3}{t(1-\frac{1}{\eta})\alpha l -|L_{\li{AB}}|(1-\frac{1}{\eta})l -(l_1 +l_2 +l_3)}\\ \leq \frac{l_1 +l_2 +l_3}{(1-\frac{1}{\eta})(\alpha -|L_{\li{AB}}|)l - (l_1 +l_2 +l_3))}.
\end{align}
Thus we have that,
\begin{align}
    r_1 \leq (r_1 +r_2)\frac{(l_1 +l_2 +l_3)}{(1-\frac{1}{\eta})(\alpha -|L_{\li{AB}}|)l}.
\end{align}
Further, as before, from the geometry of \cref{fig:third}, we have $r_1 +r_2= (1- \frac{1}{\eta})l$. 

Thus it is enough to require that,

\begin{align*}
    l-\frac{(l_1 +l_2 +l_3)}{(\alpha -|L_{\li{AB}}|)}>\frac{(l_1 +l_2 +l_3)}{\alpha}\Rightarrow l> (l_1 +l_2 +l_3)\Big( \frac{1}{(\alpha -|L_{\li{AB}}|)} +\frac{1}{\alpha} \Big)
\end{align*}

Thus as before, if we have

\begin{align*}
    1>(2\kappa +|L_\alpha|)\Big( \frac{1}{(\alpha -|L_{\li{AB}}|)} +\frac{1}{\alpha} \Big),
\end{align*}

then the previous inequality is satisfied.

Thus with the above, it is enough to require that:

\begin{align}\label{eq:eq43}
    1>(2\kappa +|L_\alpha|)\Big( \frac{1}{(\alpha -|L_{\alpha}|)} +\frac{1}{\alpha} \Big).
\end{align}
It remains also to satisfy the inequalities corresponding to \cref{new1,new2} in this situation. We note that in this case, $|CX|\leq |L_\alpha|l_1$, and then by the same argument as for \cref{new2}, it is enough to require in this situation that,
\begin{align}
    (\alpha - |L_\alpha|)l \geq 3\alpha |L_\alpha|\kappa l (\geq  3\alpha |L_\alpha|l_1)\Leftrightarrow (\alpha -|L_\alpha|)\geq 3\alpha|L_\alpha|\kappa.
\end{align}
In particular, this requires, 
\begin{align}\label{new8}
    \kappa \leq \frac{\alpha -|L_\alpha|}{3\alpha |L_\alpha|}.
\end{align}

\item[$\mathcal{z}_2:$] On the other hand, if one of the two lengths $l_{1},l_{3}$ is greater than $\kappa l$, and without loss of generality we consider that $$l_{1}>\kappa l,$$ then consider \cref{fig:inksone3}(b) where we consider a segment $\li{RP}\subset \Phi^{-n_0}(\gamma^{s}(x))$. If either $\li{RP} \cap \li{JA}$ or $\li{RP}\cap \li{BC}$ then we are done. 

If not, we have that $R$ is above $B$ in \cref{fig:inksone3}(b) and it is seen that,
\begin{align}\label{eq50new}
 l_4\geq l_{1}-|L_{\alpha}|l>(\kappa -|L_{\alpha}|)l
 \end{align}

In this case, we have two possibilities:

    \begin{figure}[H]
\centering
\includegraphics[width=1.0\textwidth]{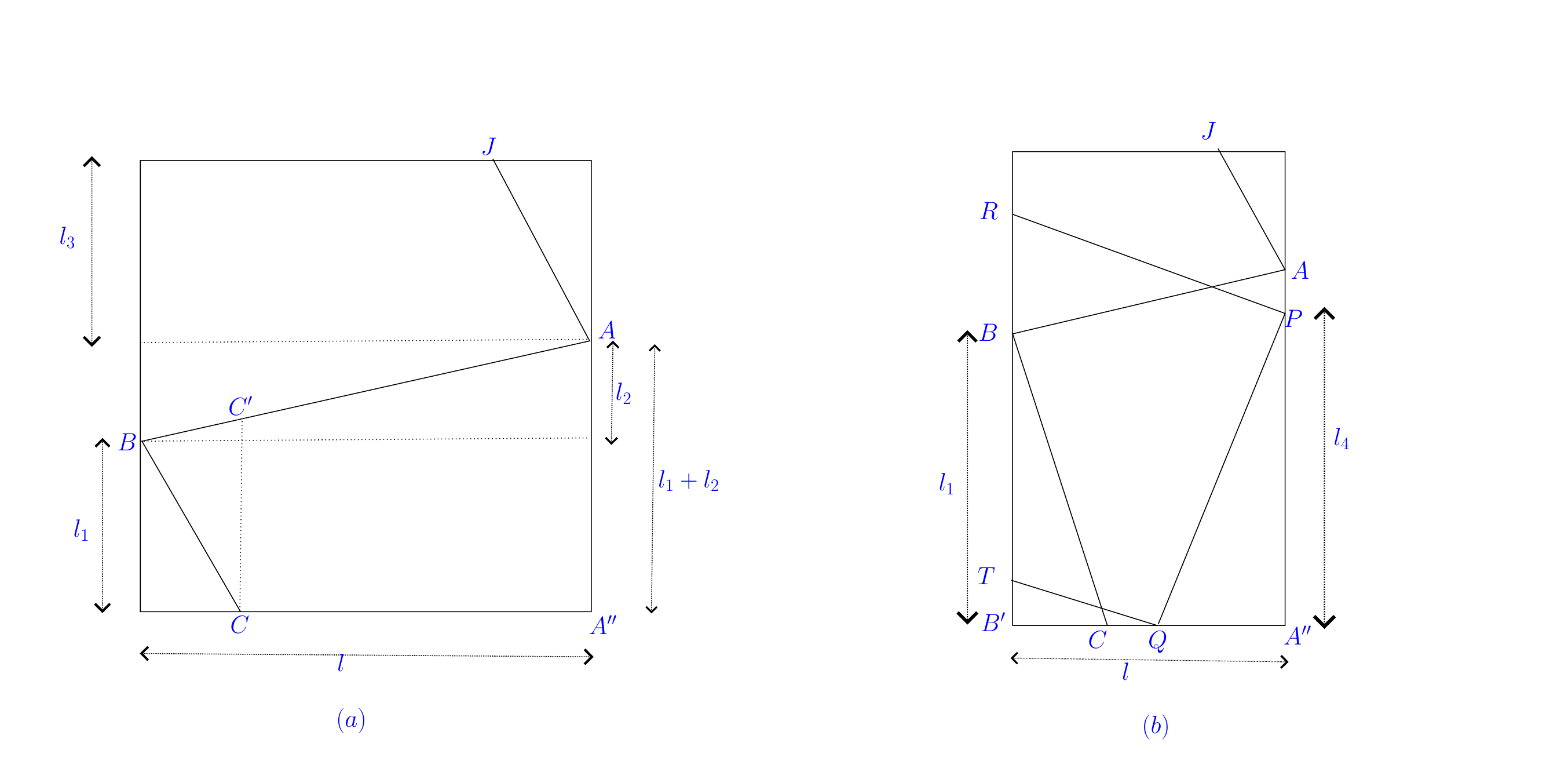}
\caption{A horizontal segment $\li{BA}\subset \FF\circ \Phi^{m_0}(\gamma^{u}(x))$ through $S$ in part (a). Part(b) shows a segment $\li{RP}\subset \Phi^{-n_0}(\gamma^{s}(y))$. Further, the segment $\li{PQ}\subset \GG^{-1}\circ\Phi^{-n_0}(\gamma^{s}(y)) $ and $\li{TQ}\subset \Phi^{-n_0 -1}(\gamma^{s}(y)). $  }
\label{fig:inksone3}
\end{figure}

\begin{enumerate}[leftmargin=*]
    \item The point $C$ lies to the left of the point $Q$ on the bottom edge of $S$. In this case, it will be enough to ensure that $\FF^{-1}(P)$ lies to the left of the left edge of $S$, and for that using \cref{eq50new} it is enough to ensure that 

    \begin{align}\label{new7}
        \alpha(\kappa-|L_{\alpha}|)>1,
    \end{align}

    in which case we can guarantee an intersection point between the segments $\li{TQ}$ and $\li{BC}$. The point of intersection is guaranteed in the case that $l_{4}<l_{1}$ since then the vertical height of $\FF^{-1}(P)$ is also less than the vertical height of $B$. If the point $P$ lies above $A$ on the right edge of $S$ then we are again done since we then have a point of intersection between $\li{RP}$ and $\li{JA}$, and so it remains to deal with the case where $P$ is below $A$ and $l_{4}\geq l_1$. In this case, with an analysis similar to that of \cref{fig:third}, we can verify (using similar triangles $\Delta(TB'Q)$ and the triangle formed by the points $Q$, $\mathbf{F}^{-1}(P)$ and the horizontal projection of $\mathbf{F}^{-1}(P)$ on the horizontal line), 
    
    $$|\li{TB'}|=\frac{l_4 |\li{B'Q}|}{\alpha l_4 -|\li{QA''}|}<\frac{l_4 l}{\alpha l_4 -l}.$$

    In this case, it is enough to ensure that, 

    \begin{align}
       |\li{TB'}|= \frac{l_4 l}{\alpha l_4 -l}< l_1 \Rightarrow \Big(\sqrt{\alpha}l_1 -\frac{l}{\sqrt{\alpha}} \Big) \Big(\sqrt{\alpha}l_4 -\frac{l}{\sqrt{\alpha}} \Big) > \frac{l^{2}}{\alpha}.
    \end{align}

    This is ensured by taking,

    \begin{align}
        \Big(\sqrt{\alpha}l_1 -\frac{l}{\sqrt{\alpha}} \Big)\geq \frac{l}{\sqrt{\alpha}}\Rightarrow l_1 > \frac{2}{\alpha}l
    \end{align}

    Thus combining the earlier estimate, it is enough to require that

    \begin{align}\label{eq:eq47}
       \kappa\geq \frac{2}{\alpha}.
    \end{align}

    \item The point $C$ lies to the right of $Q$.\footnote{This case is not shown in the figure.} For this to happen, at least one of $|\li{B'C}|$ or $|\li{QA''}|$ is greater than or equal to $l/2$. In case $|\li{QA''}|\geq l/2$, then we also have $l_4 \geq l/(2|L_{\alpha}|)$. Considering that $P$ lies below $A$, this also means that $l_1\geq l/(2|L_\alpha|)-|L_{\alpha}|l$. 

In case $|\li{B'C}| \geq l/2$, we get by the same reasoning that $l_1\geq l/(2|L_\alpha|)$ and then by hypothesis that $l_{4}\geq l/(2|L_\alpha|)-|L_{\alpha}|l$. 

In case we have $l_4\geq l_1$, we cut off the segment $\li{QP}$ up-to the vertical height of $l_1$ and subsequently call this cut off segment of height $l_1$ as $\li{QP}$ itself.

    Now again, we look at the equivalents of \cref{eq:eq32} and \cref{eq:eq36}  for the segment $\li{QP}$ which ensures that we get a horizontal segment through $S$ belonging to $\Phi^{-n_1}(\gamma^{s}(y))$ and further because of the restriction imposed in the previous paragraph, this horizontal segment must intersect the segment $\li{BC}$ and then we are done.

    For the analog of \cref{eq:eq32}, in this case it is enough to ensure that,  

    \begin{align}\label{eq:eq48}
        \alpha\Big( \frac{1}{2|L_\alpha|}-|L_\alpha| \Big)l-2l > l \Rightarrow \alpha > \frac{3}{(\frac{1}{2|L_\alpha|}-|L_\alpha|)}.
    \end{align}

    Here we are crudely bounding by $l$, from above, the distances of the points $P,Q$ from the left edge of $S$ in \cref{fig:inksone3}. \footnote{The horizontal distance of $P$ from the left edge of $S$ and the horizontal distance of $T$ from the straight line $\li{PQ}$ correspond to the quantities that appear in \cref{eq:eq32}. }

    Further, the argument in this case corresponding to the one preceding that of equation \cref{eq:eq36}, requires us to ensure:

    \begin{align*}
        \Big( \frac{1}{2|L_\alpha|}-|L_\alpha| \Big)l-\frac{l}{(\alpha -L_{\li{PQ}})}>\frac{l}{\alpha}\Rightarrow  \Big( \frac{1}{2|L_\alpha|}-|L_\alpha| \Big)>\frac{1}{\alpha} +\frac{1}{(\alpha-|L_{\li{PQ}}|)},
    \end{align*}

    where $L_{\li{PQ}}$ is the slope of the segment $PQ$ within the cone in which it lies.

    Thus because of the above, it is enough to ensure that,

    \begin{align}
        \Big( \frac{1}{2|L_\alpha|}-|L_\alpha| \Big)>\frac{1}{\alpha} +\frac{1}{(\alpha-|L_\alpha|)}.
    \end{align}
Lastly, we need to find the constraints corresponding to \cref{new1,new2} in this case. We note that $|\li{TB'}|\leq l|L_\alpha|$, and so following the argument preceding \cref{new2}, it is enough to require that, 
\begin{align}
    (\alpha -|L_\alpha|)l_4 \geq (3\alpha)l|L_{\alpha}|.
\end{align}
Thus it is enough to require that, 
\begin{align}\label{new3}
    (\alpha -|L_\alpha|)(\frac{1}{2|L_\alpha|}-|L_\alpha|)\geq 3\alpha|L_\alpha|.
\end{align}

\end{enumerate}

\end{enumerate}

    \item A vertical segment through $S$ belonging to $\Phi^{m_0}(\gamma^{u}(x))$. If now we also have a horizontal segment through $S$ belonging to $\Phi^{-n_0}(\gamma^{s}(y))$ for some integer $n_0$, then we are done, with having obtained an intersection point. Otherwise we assume that there is a vertical segment through $S$ belonging to $\GG^{-1}\circ \Phi^{-n_0}(\gamma^{s}(y))$. Now the situation is entirely analogous to the two cases $\mathcal{z}_1,\mathcal{z}_2$ before and the bounds on $\alpha$ that we get are exactly the same as those obtained from these two earlier cases $\mathcal{z}_1,\mathcal{z}_2$: we can either obtain a horizontal segment through $S$ belonging to $ \Phi^{-n_0-1}(\gamma^{s}(y))$ similar to the case $\mathcal{z}_1$, or in the remaining case find an intersection point analogous to the argument in $\mathcal{z}_2$ above.
\end{enumerate}

\subsubsection{Case b.} In the limiting case, there is a unique sequence of segments under consideration. To see this, refer to \cref{fig:inksone}. We note that under sufficiently many iterations, our four successive segments will be arbitrarily close the segments of the rectangle $\overline{ABCD}$. Without loss of generality, from now on, we assume that the four successive segments in consideration are precisely $\li{AB},\li{BC},\li{CD},\li{DA}$.

We have

\begin{align*}
    l_7+l_8 =l_3 +l_4  \\
    \Rightarrow l_6 |L_\alpha|+\frac{l_1}{|L_\alpha|} =|L_\alpha|l_2 +\frac{l_5}{|L_\alpha|} \\
    \Rightarrow l_6 |L_\alpha|^2+l_1 =|L_\alpha|^2 l_2 +l_5 \\
    \Rightarrow  l_1 -l_5 =|L_\alpha|^{2}(l_2 -l_6)\\
\end{align*}

\begin{align*}
    \text{also},\  l_1 +l_2 =l_5 +l_6 \\
    \Rightarrow  l_1 -l_5 =-(l_2 -l_6).
\end{align*}

This forces us to conclude that in \cref{fig:inksone}, 

\begin{align*}
    (|L_\alpha|^{2}+1)(l_2 -l_6)=0\implies l_2=l_6, \text{and further,}\ l_1=l_5.
\end{align*}

A similar argument lets us conclude that 

\begin{align*}
    l_8=l_4, \ l_7=l_3,
\end{align*}

and further since the ratios of $l_i/l_{i-1}=|L_\alpha|$ for $i=1,3,5,7$ (and identifying $l_0\equiv l_8$), we are forced to have a unique configuration of such segments.

Now we have either a horizontal segment through $S$, belonging to $\Phi^{-n_0}(\gamma^{s}(y))$ for some integer $n_0$, or we have a vertical segment through $S$, belonging to $\GG^{-1}\circ \Phi^{-n_0}(\gamma^{s}(y))$ for some integer $n_0$ (this latter case is not shown in \cref{fig:inksone}). One situation corresponding to the first case is shown in \cref{fig:inksone}. 

\begin{enumerate}

\item In the first case, referring to \cref{fig:inksone}, if the segment $\li{PQ}$ intersects either one of the `vertical' segments $\li{AB},\li{CD}$, then we are done. Otherwise, we have the case actually shown in \cref{fig:inksone}. In this case, we have the point $C$ is above the point $Q$ on the right edge of $S$ and the point $P$ is above $A$ on the left edge of $S$. In this case, considering the slope of $\li{PQ}$ within it's cone it is easy to see that the maximum vertical separation between the points $P,Q$ is $(l_1 +l_2)|L_\alpha|$, and this means that 

\begin{align*}
    l_8> l_3> l_8 -(l_1 +l_2)|L_\alpha|
\end{align*} 

From here we get $l_8< l_3 +(l_1 +l_2)|L_\alpha|$ and thus $l_1/|L_\alpha|< (l_1 +2l_2)|L_\alpha|$ from which we get 

\begin{align}\label{eq:another}
    l_{2}>l_{1}\frac{(1-|L_\alpha|^{2})}{2|L_\alpha|^{2}}.
\end{align}

Relative to the top edge of $S$, we look at the position of the segment $\li{BC}$. Under certain restriction on $\alpha$, we will be able to ensure a vertical segment through $S$ belonging to $\Phi^{m_1}(\gamma^{u}(x))$ for some integer $m_1$.

For the analog of \cref{eq:eq32} in this situation, it is enough to ensure (noting that $l_8=l_4$),

\begin{align*}
    \alpha l_2 -(1+|L_\alpha|^{2})l_4 -(l_3 +l_4)>(l_3 +l_4)\Rightarrow \alpha l_2 > (3+|L_\alpha|^{2})l_8 +2l_3.
\end{align*}

Also noting that $l_3 =|L_\alpha|l_2$, from the above we have the requirement of, 

\begin{align}
    (\alpha -2|L_\alpha|) l_2 > (3+|L_\alpha|^{2})l_8.
\end{align}

Along with \cref{eq:another} and the fact that $l_1=|L_\alpha|l_8$ and the above, it is enough to require that:

\begin{align}\label{eq:eq52}
    (\alpha -2|L_\alpha|)(1-|L_\alpha|^2)>2|L_\alpha|(3+|L_\alpha|^2).
\end{align}

For the analog of \cref{eq:eq36} in this situation, we have, if not a vertical segment through $S$, the requirement of 

\begin{align*}
    \alpha \Big( l_2 - \frac{l_3 +l_8}{\alpha -|L_\alpha|} \Big)>(l_3 +l_8)\Rightarrow l_2 >\Big( \frac{l_3 +l_8}{\alpha} \Big) +\Big( \frac{l_3 +l_8}{\alpha-|L_\alpha|} \Big).
\end{align*}

Noting again that $l_3=|L_\alpha|l_2$ and $l_1=|L_\alpha|l_8$, we get from the above with a bit of algebra that 

\begin{align*}
    |L_\alpha|\alpha(\alpha-|L_\alpha|)l_2>(|L_\alpha|^2 l_2 +l_1)(2\alpha -|L_\alpha|)\\ 
    \Rightarrow l_2>l_1\Bigg( \frac{2\alpha -|L_\alpha|}{\alpha^2 |L_\alpha|-3\alpha|L_\alpha|^2 +|L_\alpha|^3} \Bigg)
\end{align*}

Along with \cref{eq:another}, it is enough to require that: 

\begin{align}\label{eq:eqq53}
    \frac{(1-|L_\alpha|^{2})}{2|L_\alpha|}>\Bigg( \frac{2\alpha -|L_\alpha|}{\alpha^2 -3\alpha|L_\alpha| +|L_\alpha|^2} \Bigg)
\end{align}
Lastly, we also need to ensure the condition corresponding to \cref{new1,new2} in this situation.
We need,
\begin{align}
    (\alpha -|L_\alpha|)l_2 \geq (3\alpha)l_5=3\alpha l_1
\end{align}
Using \cref{eq:another}, it is thus enough to require that,
\begin{align}\label{new4}
    (\alpha -|L_\alpha|)\frac{(1-|L_\alpha|^{2})}{2|L_\alpha|^{2}}\geq 3\alpha.
\end{align}

\begin{figure}[H]
\centering
\includegraphics[width=0.6\textwidth]{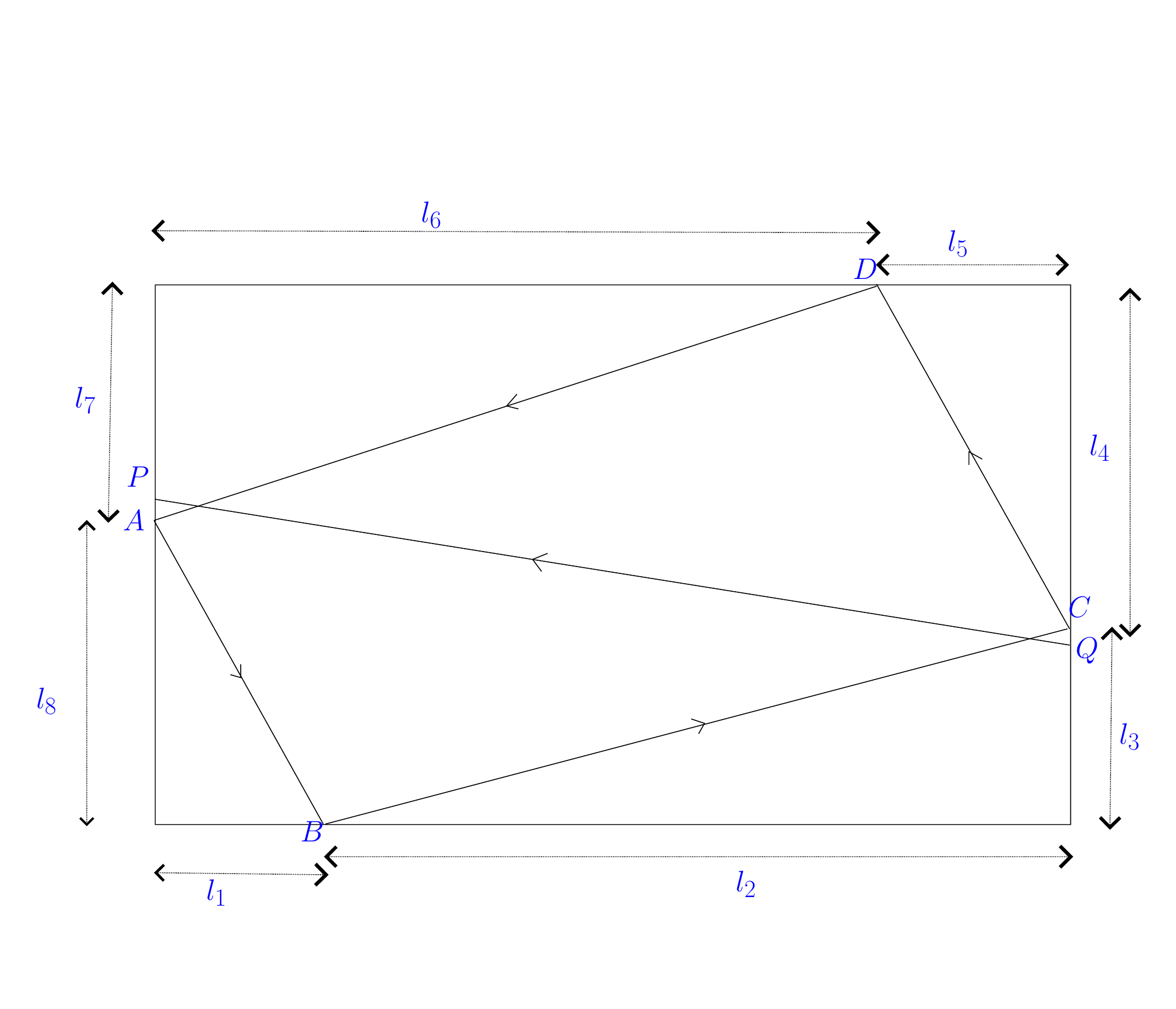}
\caption{For the forward iterates of $\gamma^{u}(x)$ we have the successive segments which approach in the limit the rectangle $\overline{ABCD}$, with the `vertical' segments $\li{AB}\subset \Phi^{m_1}(\gamma^{u}(x)), \li{CD}\subset \Phi^{m_1+1}(\gamma^{u}(x))$ and the `horizontal' segments $\li{BC}\subset \FF\circ\Phi^{m_1}(\gamma^{u}(x)),\li{AD}\subset \FF\circ\Phi^{m+1}(\gamma^{u}(x))$ for some integer $m_1$. Further, we have the segment $\li{PQ}\subset \Phi^{-n}(\gamma^{s}(y))$ for some integer $n$. }
\label{fig:inksone}
\end{figure}

\bigskip

\item Now we have remaining the case of a vertical segment $\li{PQ}$ through $S$, with $\li{PQ}\subset \GG^{-1}\circ \Phi^{-n_1}(\gamma^{s}(y))$ for some integer $n_1$. This case is shown in \cref{fig:inksone4}. If the `horizontal' segments $\li{TP}$ intersects the vertical segment $\li{AB}$, then we are done. 

Otherwise, it is enough to consider a case as shown in \cref{fig:inksone4}, with the segment $\li{PQ}\subset \GG^{-1}\circ \Phi^{-n_1}(\gamma^{s}(y))$ and with $\li{TP}\subset \Phi^{-n_1 -1}(\gamma^{s}(y))$ .

First note that by construction, in \cref{fig:inksone4}, we must have that $C$ lies below $A$, and thus also that simultaneously $U$ cannot lie below $C$ and $T$ lying above $A$. It is then apparent by the symmetry, that it is enough to consider the case where $T$ is above $A$ on the left edge of $S$ and further that the point $Q$ lies to the right of $D$.  If this is not the case, then we are either in the case where $\li{QU}$ intersects $\li{DC}$ where we are done, or we have the case where $U$ lies below $C$ on the right edge of the square. In that case, looking at \cref{fig:inksone4}, one places the segment $\li{PQ}$ in place of $\li{SR}$ in which case if $P$ (in place of $R$) lies to the right of $B$ on the bottom edge, then also we're done since $\li{PT}$ would then intersect $\li{AB}$. So the only case that remains to consider is where $Q$ lies to the right of $D$ on the top edge, or equivalently by symmetry if $R$ lies to the left of $B$ on the bottom edge.

In this case we will construct below a segment $\li{RS}\subset \GG^{-1}\circ \Phi^{-n_2}(\gamma^{s}(y)) $, for some integer $n_2>n_1$, and where $S$ on the top edge lies to the left of the point $D$. This will either force a point of intersection between the segments $\li{SZ}$ and $\li{DC}$ or $Z$ lies below $C$ on the right edge of $S$, and then further onwards we can get a segment $\li{EF}\subset \Phi^{m_2}(\gamma^{u}(x))$ which will intersect at least one of $\li{TP}$ or $\li{SZ}$ and we would be done.

\begin{figure}[h]
\centering
\includegraphics[width=0.65\textwidth]{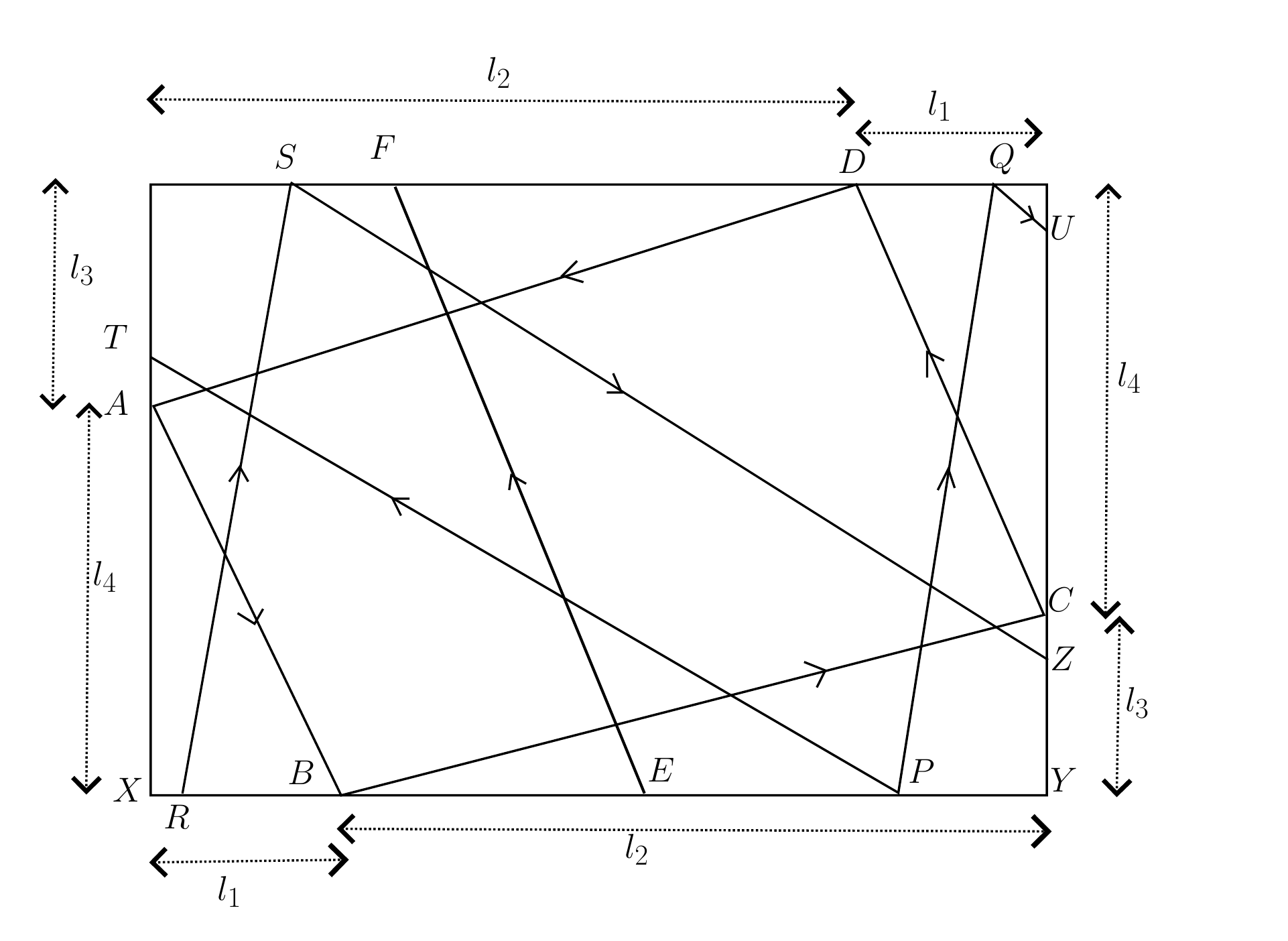}
\caption{The segments constituting the rectangle $\overline{ABCD}$ are as in \cref{fig:inksone}. Now we are considering a segment $\li{PQ}\subset \GG^{-1}\circ \Phi^{-n_1}(\gamma^{s}(y))$ for some positive integer $n_1$. We further have $\li{TP}\subset \Phi^{-n_1 -1}(\gamma^{s}(y))$ and then we get a segment $\li{RS}\subset \GG^{-1}\circ \Phi^{-n_2}(\gamma^{s}(y)) $ (employing the argument preceding \cref{eq:eq32}) for some integer $n_2>n_1$, in such a way that the segment $\li{SZ}\subset \Phi^{-n_2-1}(\gamma^{s}(y))$ either intersects $\li{CD} \subset  \Phi^{m_1+1}(\gamma^{u}(x))$ or else we will find a vertical segment through $\li{EF}\subset \Phi^{m_2}(\gamma^{u}(x))$ through $S$ which is now guaranteed to intersect at least one of the segments $\li{TP},\li{SZ}$ and we would be done. }
\label{fig:inksone4}
\end{figure}

Since $T$ lies above $A$ on the left edge of $S$ in \cref{fig:inksone4}, we have, with $L_{\li{PQ}}$ being the slope of $\li{PQ}$ within it's cone, that $l_4(\alpha-|L_\alpha|)\leq |\li{TX}|(\alpha+ L_{\li{PQ}})=|\li{PX}|\leq (l_1 +l_2)$, and thus also noting that $l_1 =|L_\alpha| l_4$, we have

\begin{align}\label{eq:eq53}
    l_4(\alpha-2|L_\alpha|)\leq l_2.
\end{align}

Now using the segment $\li{TP}$ relative to the top edge of $S$, we employ a modified argument preceding \cref{eq:eq32} in order to ensure that we have the requisite vertical segment $\li{RS}$ through the square $S$. First note that clearly 

\begin{align}\label{eq:eq67}
|\li{PX}|\geq l_1 +l_2 -(l_1 + |L_\alpha (l_3 +l_4)| = (l_2- |L_\alpha|(l_3 +l_4)).
\end{align}

In this case, using the argument prior to \cref{eq:eq32}, it is enough to require,

\begin{align*}
    \alpha(l_2- |L_\alpha|(l_3 +l_4))-(l_3 +l_4)-l_3 (1+|L_\alpha|^2)> (l_3 +l_4) \\
    \Rightarrow \alpha(l_2- |L_\alpha|(l_3 +l_4))>l_3 (3+|L_\alpha|^2)+2l_4.
\end{align*}

Also noting that, $l_3=|L_\alpha|l_2$, the above reduces to requiring

\begin{align}
    l_2 \Big(\alpha(1-|L_\alpha|^2) -|L_\alpha|(3+|L_\alpha|^2)>l_4 (2+\alpha|L_\alpha|) \Big)
\end{align}

Thus combining \cref{eq:eq53}, it is enough to require that 

\begin{align}\label{eq:eq56}
    (\alpha-2|L_\alpha|)>\frac{(2+\alpha|L_\alpha|)}{\Big(\alpha(1-|L_\alpha|^2) -|L_\alpha|(3+|L_\alpha|^2)\Big)}.
\end{align}

Further, using an argument similar to the one preceding \cref{eq:eq36} in this case, it is easily seen that we will be done with a further requirement that:

\begin{align*}
    \alpha\Big( |\li{PX}|-\frac{(l_3 +l_4)}{(\alpha-|L_\alpha|)} \Big)>(l_3 +l_4)
\end{align*}

and thus it is enough to require that:

\begin{align*}
    (l_2 -|L_\alpha|(l_3 +l_4))>\frac{(l_3 +l_4)}{(\alpha-|L_\alpha|)}+\frac{(l_3 +l_4)}{\alpha}\\ \Rightarrow l_2\Big( 1- \frac{|L_\alpha|}{(\alpha -|L_\alpha|)}-\frac{|L_\alpha|}{\alpha}\Big) >l_4 \Big( \frac{1}{\alpha}+\frac{1}{\alpha-|L_\alpha|}+|L_\alpha|  \Big).
\end{align*}

Thus combining with \cref{eq:eq53}, it will be enough to ensure that:

\begin{align}\label{eq:eq57}
    (\alpha-2|L_\alpha|)>\frac{(\alpha-|L_\alpha|)+\alpha +|L_\alpha|\alpha(\alpha-|L_\alpha|)}{\alpha(\alpha-|L_\alpha|)-\alpha|L_\alpha|-(\alpha-|L_\alpha|)|L_\alpha|}.
\end{align}

Further, we use the argument of \cref{new1,new2} in this situation, and get noting that $T$ lies above $A$ on the left edge,
\begin{align}
    (\alpha -|L_\alpha|)|\li{PX}|\geq 3\alpha l_3 |L_\alpha|.
\end{align}
In this case, using \cref{eq:eq67}, it is enough to require that,
\begin{align}
    (\alpha -|L_\alpha|)(l_2 -|L_\alpha|(l_3 +l_4))\geq 3\alpha l_3 |L_\alpha|.
\end{align}
Using that, $l_3 =|L_\alpha| l_2$,  ,it is enough to require that
\begin{align}
    (\alpha -|L_\alpha|)(l_2 -|L_\alpha|(|L_\alpha |l_2 +l_4))\geq 3\alpha l_2|L_\alpha|^{2}\\ \Leftrightarrow (\alpha -|L_\alpha|) l_2 (1-|L_\alpha|^{2}) \geq 3\alpha l_2 |L_\alpha|^{2} +l_4 |L_\alpha|(\alpha -|L_\alpha|).
\end{align}
Now noting that and $l_2 \geq l_4 (\alpha -2|L_\alpha|)$, it is enough to require that,
\begin{align}\label{new5}
    (\alpha -|L_\alpha|) (1-|L_\alpha|^{2}) \geq 3\alpha |L_\alpha|^{2} + \frac{|L_\alpha|(\alpha -|L_\alpha|)}{(\alpha -2|L_\alpha|)}.
\end{align}
We note that as a result we have the vertical segment $\li{RS}$, and that the point $S$ lies to the left of $D$ on the top edge of the square. By looking at the position of the segment $\li{BC}$ relative to the top edge, with arguments similar to those preceding \cref{eq:eq32,eq:eq36}, we ensure for strong enough twists that there is the vertical segment $\li{EF}\subset \Phi^{m_2}(\gamma^{u})(x)$ for some integer $m_2 >m_1$.

With an analog of the argument preceding \cref{eq:eq32}, it is enough to require that:

\begin{align*}
    \alpha l_2 -(l_3 +l_4)-l_4 (1+|L_\alpha|^2)>(l_3 +l_4),\\
   \Leftrightarrow \alpha l_2>2(l_3 +l_4) +l_4 (1+|L_\alpha|^2),\\
    \Leftrightarrow (\alpha-2|L_\alpha|)l_2 >l_4 (3+|L_\alpha|^2),
\end{align*}

and thus with \cref{eq:eq53}, it is enough to require that

\begin{align}\label{eq:eq58}
    (\alpha-2|L_\alpha|)^2> (3+|L_\alpha|^2)
\end{align}

Further, for an analog of the argument preceding \cref{eq:eq36}, we either have the vertical segment through the square in which case we are done, otherwise it is enough to have

\begin{align*}
    \Big(l_2 - \frac{(l_3 +l_4)}{(\alpha-|L_\alpha|)}\Big)>\frac{(l_3 +l_4)}{\alpha}\\
    \Rightarrow l_2 \big( \alpha^2 -3\alpha|L_\alpha|+|L_\alpha|^2 \big)>l_4\big( 2\alpha -|L_\alpha| \big)
\end{align*}

Along with \cref{eq:eq53}, it is enough to ensure that 

\begin{align}\label{eq:eq59}
    (\alpha-2|L_\alpha|)>\frac{(2\alpha -|L_\alpha|)}{(\alpha^{2}-3\alpha|L_\alpha|+|L_\alpha|^2)}.
\end{align}
We now need the argument corresponding to \cref{new1,new2}. We need,
\begin{align}
     (\alpha -|L_\alpha|)l_2\geq 3\alpha l_1 =3\alpha l_4 |L_\alpha|.
\end{align}
Again, using the fact that $l_2\geq l_4 (\alpha -2|L_\alpha|)$, we see that it is enough to require that, 
\begin{align}\label{new6}
    (\alpha -|L_\alpha|)(\alpha-2|L_\alpha|)\geq 3\alpha |L_\alpha|.
\end{align}
Lastly we need to ensure that the segment $\li{EF}$ does intersect at least one of the segments $\li{TP}$ or $\li{SZ}$; for which it is enough to ensure that the magnitude of the slope of the line segment $\li{SP}$ is greater than $|L_\alpha|$.

In this case, we note that the horizontal separation of the points $S,P$ is clearly bounded from below by $(l_1 +l_2)-2(l_1 +|L_\alpha|(l_3 +l_4))$ and thus it is enough to require that 

\begin{align*}
    \frac{l_2 -l_1 -2|L_\alpha|(l_3 +l_4)}{(l_3 +l_4)}\geq |L_\alpha|\\ \Rightarrow l_2 -l_1 \geq 3|L_\alpha|(l_3 +l_4),
\end{align*}

and noting that $l_1=|L_\alpha|l_4$ and $l_3=|L_\alpha|l_2$, the above gives

\begin{align*}
    l_2 (1+3|L_\alpha|^2)>4l_4 |L_\alpha|.
\end{align*}

Combining with \cref{eq:eq53}, it is enough to require that 

\begin{align}\label{eq:eq60}
    (\alpha-2|L_\alpha|)(1+3|L_\alpha|^2)>4|L_\alpha|.
\end{align}
 
\end{enumerate}

\subsubsection{Case c.} This case is entirely analogous to the previous Case b and all the bounds on $\alpha$ that we get here would be exactly the same as in the Case b.

\subsubsection{Case d.} 
\begin{figure}[h]
\centering
\includegraphics[width=0.8\textwidth]{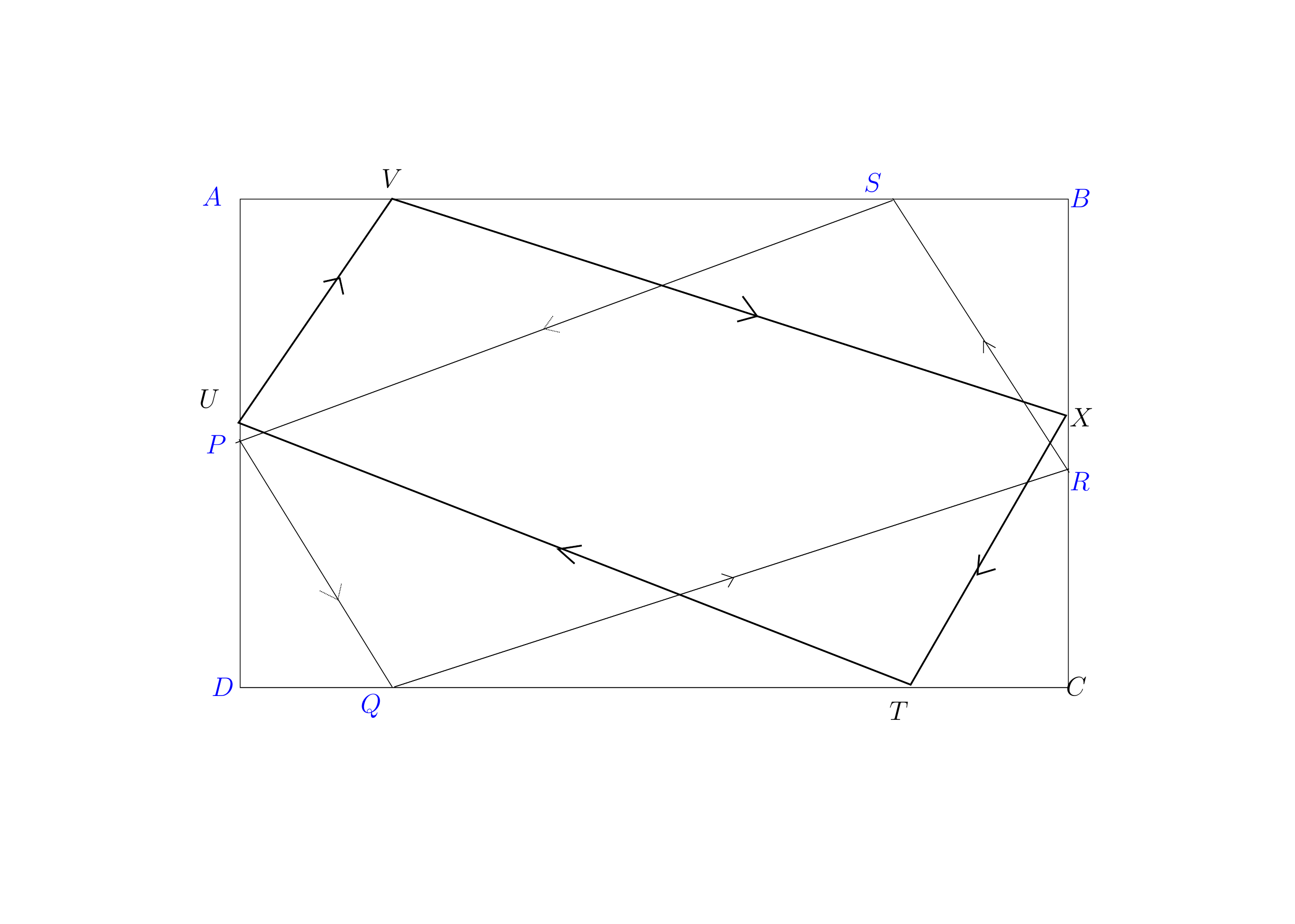}
\caption{For the forward iterates of $\gamma^{u}(x)$ we have the successive segments which approach in the limit the rectangle $\overline{ABCD}$, with the `vertical' segments $\li{AB}\subset \Phi^{m}(\gamma^{u}(x)), \li{CD}\subset \Phi^{m+1}(\gamma^{u}(x))$ and the `horizontal' segments $\li{BC}\subset \FF\circ\Phi^{m}(\gamma^{u}(x)),\li{AD}\subset \FF\circ\Phi^{m+1}(\gamma^{u}(x))$ for some integer $m$. The rectangle $\overline{XTUV}$ in the figure represents the corresponding set of sequences for the backward iterates of $\gamma^{s}(y)$. }
\label{fig:inksone2}
\end{figure}

In this case, we have four successive segments within $S$ belonging to forward iterates of $\gamma^{u}(x)$ and four successive segments  within $S$ belonging to backward iterates of $\gamma^{s}(y)$. Further, these two limiting rectangles are reflections of each other about a vertical line through the middle of $S$. 

We refer to \cref{fig:inksone2}. In this case, as noted earlier, we have $|\li{RC}|\leq |\li{PD}|$. In case we have strictly $|\li{RC}|< |\li{PD}|$, then we also necessarily have $|\li{UD}|<|\li{PD}|$ and that $U$ is actually below $P$, when the rectangle $\overline{XTUV}$ is taken to be the limiting set of four segments for the backward iterates of the stable manifold $\gamma^{s}(y)$. Thus after sufficiently many iterates, when the segments belonging to the forward iterates of the unstable manifold $\gamma^{u}(x)$ and the segments belonging to the backward iterates of the stable manifold $\gamma^{s}(y)$ are sufficiently close respectively to the rectangles $\overline{ABCD}$ and $\overline{XTUV}$ then we are forced to get requisite points of intersection between the segments $\li{PQ},\li{TU}$ and also the segments $\li{VX},\li{SR}$.

The only case that remains is when $|\li{RC}|= |\li{PD}|$. This can in any case correspond to a Lebesgue measure zero set of points $\{(x,y)\in S^2\}$ and thus makes no difference to the argument. Even then we can guarantee a requisite point of intersection in this case. In this case, we assume without loss of generality that the four forward iterates of  $\gamma^{u}(x)$ are arbitrarily close to the segments of the rectangle $\overline{PQRS}$. Further suppose, as in \cref{fig:inksone2}, that the segment $\li{XT}$ belongs to some backward iterate of  $\gamma^{u}(x)$ and that further, the point $X$ lies above $R$ and then also the point $U$ lies above $P$, in which case we will be then forced to have a point of intersection between $\li{SR}$ and $\li{VX}$ and we would be done. Note that in this instance, the set of segments of $\overline{XTUV}$ as shown in \cref{fig:inksone2} do not constitute the limiting rectangle. In fact, this may actually not even be a proper rectangle and the image of $U$ under $\FF^{-1}$ will not in general coincide with $X$ but we will have $|\li{BX}|\leq |\li{AU}|$. If the point $X$ lies above $R$ and the point $U$ lies below $P$, then we have a requisite point of intersection between the segments $\li{TU}$ and $\li{PQ}$ and we are done.

In the remaining case, if $X$ lies below $R$ then again we are forced to have a point of intersection between $\li{TU}$ and $\li{PQ}$ and we are also done.

\bigskip

\subsection{Optimization in Case(ii).}

In the end we have to optimize the parameters $\beta_1,\beta_2$, and also take into account the several constraints on $\alpha$, arising from \cref{eq:eq13,eq:eq16,eq:eq18,eq:eq22,eq:eq23,eq:eq32,eq:eq36} from Section 2.3(Part 2), along with \cref{eq:eq42,eq:eq43,,eq:eq47,eq:eq48,eq:eq52,eq:eqq53,eq:eq56,eq:eq57,eq:eq58,eq:eq59,eq:eq60,new1,new2,new3,new4,new5,new6,new7,new8}  from Section 2.4(Part 3). We require the paramters $\beta_1,\beta_2$ to be the minimum possible so the there exists a $\delta>1$ so that the inequalities from Section 2.3 hold, and also ensure the twist parameter satisfies the constraints imposed by the further inequalities of Section 2.4. 

\begin{enumerate}[leftmargin=*]

\item We first deal with the constraints of Section 2.3.

From \cref{eq:eq22}, we have the requirement:

\begin{align}\label{req3}
    \beta_2 > \frac{\delta}{L_5+\alpha}+  \frac{\delta}{L_7+\alpha}.
\end{align}

From \cref{eq:eq36}, we have the requirement:

\begin{align}\label{req1}
    \beta_2 > \frac{\delta}{\alpha-|L_5|}+  \frac{\delta}{\alpha}.
\end{align}

From \cref{eq:eq32}, we have the requirement:

\begin{align}\label{req2}
    \beta_2 > \frac{\delta}{\alpha -|L_5|}\Big(  1+\frac{2+|L_5 L_6|}{(L_6 +\alpha)} \Big).
\end{align}
Also recall from \cref{new2}, the constraint:
\begin{align}\label{req55}
    \beta_2 \geq \frac{3\alpha |L_\alpha|\delta}{(\alpha- |L_\alpha|)^{2}}.
\end{align}

From \cref{eq:eq13}, we have the requirement:

\begin{align}\label{req3}
     \beta_1 >\frac{\delta}{L_2 +\alpha}.
\end{align}

From \cref{eq:eq16}, we have the requirement:

\begin{align}\label{req4}
    \beta_1 >\frac{\delta}{L_3+\alpha} +\frac{\delta}{L_2 +\alpha}.
\end{align}

From \cref{eq:eq18}, we have the following constraint on $\alpha$:

\begin{align}\label{req5}
      (L_3+\alpha)>\delta\Big(1+\frac{L_3 +\alpha}{L_1+\alpha}\Big)
\end{align}

For the above to hold, for some $\delta>1$, note that it is enough to have $(L_3 +\alpha)(L_1 +\alpha)> L_1 +L_3 +2\alpha$, and since we have $(L_3 +\alpha)(L_1 +\alpha)>(L_\alpha +\alpha)^{2}>2\alpha> L_1 +L_3 +2\alpha$, it is enough to require that:

\begin{align}\label{eqp}
   (L_\alpha +\alpha)^{2}> 2\alpha.
\end{align}

We note that \cref{eqp} is true when $\alpha>2.783$.

We note that the constraint from \cref{eq:eq23} is the same as the one above from \cref{eq:eq18}.

Noting that we have uniformly $L_\alpha\leq L_i \leq 0$, for all $i=1,\dots,9$, to require \cref{req1,req2,req55} it is enough to require that:

\begin{align}\label{eqq}
    \beta_{2} \geq \text{max}\Big(\frac{\delta}{\alpha-|L_\alpha|}+\frac{\delta}{\alpha}, \frac{\delta}{\alpha-|L_\alpha|}\Big(1+ \frac{2+L_{\alpha}^{2}}{\alpha-|L_\alpha|} \Big),\frac{2\delta}{\alpha-|L_\alpha|},  \frac{3\alpha |L_\alpha|\delta}{(\alpha- |L_\alpha|)^{2}}\Big)\\ =  \text{max}\Big( \frac{\delta}{\alpha-|L_\alpha|}\Big(1+ \frac{2+L_{\alpha}^{2}}{\alpha-|L_\alpha|} \Big),\frac{2\delta}{\alpha-|L_\alpha|},  \frac{3\alpha |L_\alpha|\delta}{(\alpha- |L_\alpha|)^{2}}\Big)
\end{align}

On the right hand side above, we have accordingly as $\alpha\gtrless 2.66$, that:

\begin{align*}
    \frac{2\delta}{\alpha-|L_\alpha|}\gtrless \frac{\delta}{\alpha-|L_\alpha|}\Big(1+ \frac{2+L_{\alpha}^{2}}{\alpha-|L_\alpha|}\Big) \ \  \Big( \Rightarrow \alpha-|L_\alpha| \gtrless (2+L_{\alpha}^{2})\Big).
\end{align*}

Further, we also have, accordingly as $\alpha\gtrless 2.43$, that,
\begin{align}
    \frac{2\delta}{\alpha-|L_\alpha|}\gtrless   \frac{3\alpha |L_\alpha|\delta}{(\alpha- |L_\alpha|)^{2}}
\end{align}
Also, to ensure \cref{req3,req4}, it is enough to require that:

\begin{align}
    \beta_1\geq \text{max}\Big( \frac{\delta}{\alpha-|L_\alpha|},\frac{2\delta}{\alpha-|L_\alpha|} \Big) = \frac{2\delta}{\alpha-|L_\alpha|} 
\end{align}

First we restrict to $\alpha\geq 2.66$, and take the parameters $\beta_1$ and $\beta_2$ to be the minimum possible allowable values; thus, we take:

\begin{align}
    \beta_2=\beta_1=\frac{2\delta}{\alpha-|L_\alpha|},
\end{align}

Without loss of generality, can take the $\delta=1$ above for the purpose of this calculation,
\begin{align}\label{eqbeta}
    \beta_2=\beta_1=\frac{2}{\alpha-|L_\alpha|} ,
\end{align}

Note that whenever $\alpha> 2.783$, we can alter the $\delta>1$ as necessary to get that  $(1-\delta\cdot l_{v}(\gamma))>0$, since in the most general case we can verify that $l_{v}(\gamma)<\frac{2}{(L_\alpha +\alpha)}<1$, since $2<L_\alpha +\alpha$ when $\alpha>2.783$, and we restrict to this range now. As noted earlier, $L_\alpha=-(\alpha/2)+\sqrt{(\alpha/2)^{2}-1}$. \footnote{Note, as in \cite{Prz}, that the condition $2/(L_\alpha +\alpha)<1$ follows from requiring that we are not already done with having a horizontal segment through $S$}

Thus following \cite{Prz}, it is enough to require:

\begin{align}
    l_v (\gamma)=\sum\limits_{i=1}^{4} l_v (I_i)> l_v (\gamma)\Bigg(\frac{\beta_2}{\alpha(1-l_v (\gamma))}+\frac{2}{2\alpha +L_\alpha}+\frac{\beta_1}{\alpha +L_\alpha}\Bigg),
\end{align}

and thus it is enough to have,

\begin{align}\label{eqmostimp}
    1>\Big(\frac{\beta_2}{\alpha(1-\frac{2}{\alpha+L_{\alpha}})}+\frac{2}{2\alpha +L_\alpha}+\frac{\beta_1}{\alpha +L_\alpha}\Bigg),
\end{align}

with $\beta_1, \beta_2$ being given by \cref{eqbeta}. The optimal parameter is $\alpha_0=3.47$ in this case, consistent with our earlier assumption that $\alpha>2.783$ to begin with.

\item Next we verify that the constraints arising in Section 2.4 are all satisfied for the optimal twist parameters $\alpha>\alpha_0$. 

For the constraints \cref{eq:eq42,eq:eq43,new7,new8,eq:eq47} involving $\kappa$, we a-priori use an estimate of $\kappa=2/3$, and verify that the corresponding inequalities are true for $\alpha>\alpha_0 =3.47$. We verify that \cref{eq:eq42} is satisfied for $\alpha>2.69$, \cref{eq:eq43} is satisfied for $\alpha>3.46$, \cref{eq:eq47} is satisfied for $\alpha>3$.

For the remaining constraints, we check that \cref{eq:eq48} is satisfied for $\alpha>3.07$, \cref{eq:eq52} is satisfied for $\alpha>3.20$, \cref{eq:eqq53,eq:eq56} are always satisfied, \cref{eq:eq57} is satisfied for $\alpha>2.75$, \cref{eq:eq58} is satisfied for $\alpha>3.33$, \cref{eq:eq59} is satisfied for $\alpha>2.54$, and \cref{eq:eq60} is satisfied for $\alpha>2.31$.\footnote{All these inequalities were verified with Mathematica 13.1}

\end{enumerate}

We further note that the inequality corresponding to \cref{new3} is satisfied for $\alpha\geq 3.28$, the inequality corresponding to \cref{new4} is satisfied for $\alpha\geq 3.25$,  \cref{new5} is satisfied for $\alpha\geq 2.81$, and \cref{new6} is satisfied for $\alpha\geq 2.61$. 

Thus the optimal parameter is $\alpha=\alpha_0 =3.47$.



\subsection{Case (iv)}

We now deal with the simpler Case(iv) outlined in Section 2.1 for the first return having two components each intersecting $S$, as depicted by \cref{fig:fifth}(b). For this, either of the following three is enough:

\begin{enumerate}
    \item $l_{h}(I_1)\geq \delta \beta_1 l_{v}(\gamma)$,

    \item $l_{h}(I_1)\geq \delta \beta_1 l_{v}(\gamma)$,

    \item $l_h\FF(I_2)-l_h(I_2)\geq \delta \cdot l_v(\gamma)$.
\end{enumerate}

Here the constant $\beta_1$ from the analysis of Case(ii) is used in the first two equations since these situations can be easily seen to be analogous to the corresponding analysis in Case(ii). The bounds obtained in Case (iv) are better than those obtained from Case (ii), the critical twist is determined by the best possible improvement in Case(ii), and we don't work on improving the third equation above, and only require: $l_h\FF(I_2)-l_h(I_2)\geq \delta \cdot l_v(\gamma)$. 

In this case, following Equation (13) in \cite{Prz} and the earlier arguments, it is enough to require that:

\begin{align}
    1> \frac{2\beta_1}{\alpha + L_\alpha} + \frac{1}{\alpha}.
\end{align}

The optimal value from the above is $\alpha=2.95$.

Thus combining Case(ii) (and hence also Case(iii)) along with Case(iv), we get that in the most general case, $\alpha=3.47$ is the optimal twist parameter.

\bigskip

\section{Conclusion:}

 All the extensions to the arguments of \cite{Prz} here involve segments that touch one end of the central square $S$.

The methods used here can be extended to deal with usual modifications of the linked twist map, such as those discussed in Section 2 of \cite{Prz}, or in \cite{Sp}; cases with more than one linked region. Further, one might make an assumption of the twisting regions being small, in which case one can assume that $\l_{v}(\gamma)$ is negligible, and thus alter Eq. (73) for getting the optimal parameter from Case(ii). 

Further, in the analysis of determining the lower bound on the $\beta_2$ parameter, we can make adjustments in several places if one assumes that the regions $H\setminus S$ and $V\setminus S$ are large compared with $S$ itself, or imposes certain other restrictions for the parameters $D_1,D_2$ in Figure 6, this would improve the lower bound for $\beta_2$.

Also, if one works with boundary identifications where the top edge of the unit square is identified with the left or the right edge, as happens in \cite{Pat}, then again these methods can be suitably altered to improve the optimal twist parameter for which ergodicity is achieved. 

\section{Acknowledgements:} The author is thankful to Feliks Przytycki for useful feedback on this question and to Manu Mannattil for help with setting up some of the diagrams.

\section{Declarations} 
 
\subsection{Ethical Approval} 
Not applicable
 
\subsection{Competing interests }
Not applicable
 
\subsection{Authors' contributions }
Not applicable 
 
\subsection{Funding} 
Not applicable
 
\subsection{Availability of data and materials} 
Not applicable

\FloatBarrier

\end{document}